\documentclass[12pt, a4paper]{article}

\usepackage{fancyhdr}
\usepackage{epsfig}
\usepackage{cite}
\usepackage{amsmath}
\usepackage{theorem}
\usepackage{graphicx}
\usepackage{amssymb}
\usepackage{latexsym}
\usepackage{epic}
\usepackage{amsfonts,amssymb,amsopn}
\usepackage[all,2cell,dvips]{xy} \UseAllTwocells \SilentMatrices

\newcommand{\C}{\mathbb{C}}
\newcommand{\Z}{\mathbb{Z}}
\newcommand{\PP}{\mathbb{P}}

\newcommand{\Homom}{\mathrm{Hom}}

\newcommand{\Endom}{\mathrm{End}}
\newcommand{\Ker}{\mathrm{Ker}}
\newcommand{\Image}{\mathrm{Im}}
\newcommand{\Coker}{\mathrm{Coker}}
\newcommand{\Iden}{\mathrm{Id}}

\newcommand{\Spn}{\mathrm{Sp}_{n}\C}

\newcommand{\Spt}{\mathrm{Sp}_{2}\C}
\newcommand{\diag}{\mathrm{diag}}
\newcommand{\tr}{\mathrm{tr}}

\newcommand{\degree}{\mathrm{deg}}

\newcommand{\Sym}{\mathrm{Sym}}

\newcommand{\urd}{\mathcal{U}(r,d)}
\newcommand{\surL}{\mathcal{SU}(r,L)}

\newcommand{\Pic}{\mathrm{Pic}}

\newcommand{\Mt}{\mathcal{M}_{2}}
\newcommand{\Mn}{\mathcal{M}_{n}}
\newcommand{\Mxn}{\mathcal{M}(\Spn)}
\newcommand{\Mxs}{\mathcal{M}(\Spt)}

\newcommand{\Ox}{O_{X}}
\newcommand{\Kx}{K_{X}}
\newcommand{\strs}{\mathcal{O}_{X}}
\newcommand{\cans}{\mathcal{K}_{X}}

\newcommand{\ev}{\mathrm{ev}}
\newcommand{\Supp}{\mathrm{Supp}}
\newcommand{\Bs}{\mathrm{Bs}}

\newcommand{\qed}{{\ifhmode\unskip\nobreak\fi\quad\ensuremath\square}}
\newcommand{\pfe}{\PP^{5}_{E}}
\newcommand{\pfed}{( \pfe )^{*}}
\newcommand{\J}{J^{1}}
\newcommand{\je}{j_{E}}

\newcommand{\De}{D_{E}}
\newcommand{\Fix}{\mathrm{Fix}}
\newcommand{\ftp}{|4\Theta|_{+}}
\newcommand{\hJ}{\widehat{J}}
\newcommand{\Oj}{O_{J}}
\newcommand{\strsj}{{\mathcal O}_{J}}
\newcommand{\tmo}{\widetilde{(-1_{J})}}

\newcommand{\Htt}{{\mathcal G}(2\Theta)}

{\theorembodyfont{\slshape}\newtheorem{prop}{{\textbf Proposition}}}
{\theorembodyfont{\slshape}\newtheorem{thm}[prop]{{\textbf Theorem}}}
{\theorembodyfont{\slshape}\newtheorem{cor}[prop]{{\textbf Corollary}}}
{\theorembodyfont{\slshape}\newtheorem{lemma}[prop]{{\textbf Lemma}}}
{\theorembodyfont{\slshape}}
{\theorembodyfont{\slshape}}

\title{Rank four symplectic bundles without theta divisors over a curve of genus two}

\author{George H.\ Hitching}

\begin{document}

\maketitle

\begin{abstract}
The moduli space $\Mt$ of rank four semistable symplectic vector bundles over a curve $X$ of genus two is an irreducible projective variety of dimension ten. Its Picard group is generated by the determinantal line bundle $\Xi$. The base locus of the linear system $|\Xi|$ consists of precisely those bundles without theta divisors, that is, admitting nonzero maps from every line bundle of degree $-1$ over $X$. We show that this base locus consists of six distinct points, which are in canonical bijection with the Weierstrass points of the curve. We relate our construction of these bundles to another of Raynaud and Beauville using Fourier--Mukai transforms. As an application, we prove that the map sending a symplectic vector bundle to its theta divisor is a surjective map from $\Mt$ to the space of even $4\Theta$ divisors on the Jacobian variety of the curve.
\end{abstract}

%\tableofcontents

\section{Introduction}

In this section we introduce the objects we will be studying, and give a summary of the paper. Let $X$ be a complex projective curve which is smooth and irreducible of genus $g \geq 2$.\\
\\
\textbf{Definition:} A vector bundle $W \to X$ is \textsl{symplectic} (resp., \textsl{orthogonal}) if there is a bilinear nondegenerate antisymmetric (resp., symmetric) form $\omega$ on $W \times W$ with values in a line bundle $L$, which for us will often be the trivial bundle $\Ox$. Note that a symplectic bundle is necessarily of even rank.\\
\par
Let $r$ and $d$ be integers with $r \geq 1$. The moduli space of semistable vector bundles of rank $r$ and degree $d$ over $X$ is denoted $\urd$. For any line bundle $L \to X$ of degree $d$, the closed subvariety of $\urd$ of bundles with determinant $L$ is denoted $\surL$. References for these objects include Seshadri \cite{Ses1982} and Le Potier \cite{LeP1997}. The variety ${\mathcal U}(1,d)$ is the $d$th Jacobian variety of $X$, and will be denoted $J^d$; see for example Birkenhake--Lange \cite{BL2004}, Chap.\ 11, for details. The main object of interest for us is the moduli space of semistable symplectic vector bundles of rank $2n$ over $X$, which is denoted $\Mn$. We have

\begin{thm} $\Mn$ is canonically isomorphic to $\Mxn$, the moduli space of semistable principal $\Spn$-bundles over $X$, and the natural map $\Mn \to \mathcal{SU}(2n, O_{X})$ is injective. \label{Mn=Mxn}
\end{thm}
\textbf{Proof}\\
See \cite{Hit2005}, chap.\ 1, or \cite{Hit2005b} for a sketch.\\
\\
This allows us to use Ramanathan's results in \cite{Rthn1996a} and \cite{Rthn1996b} (especially Theorem 5.9) to give information about $\Mn$. We find that $\Mn$ is an irreducible, normal, projective variety of dimension $n(2n+1)(g-1)$.

\subsection*{Generalised theta divisors and the theta map}

For a semistable vector bundle $W \to X$ of rank $r$ and trivial determinant, we consider the set
\[ S(W) := \left\{ L \in J^{g-1} : h^{0}(X, L \otimes W) > 0 \right\}. \]
If $S(W) \neq J^{g-1}$ then it is the support of a divisor $D(W)$ on $J^{g-1}$ linearly equivalent to $r \Theta$, called the \textsl{theta divisor} of $W$. It is not hard to show that this only depends on the $S$-equivalence class of $W$ (see for example \cite{Hit2005}, chap.\ 6). The association $D \mapsto D(W)$ defines a rational map $\mathcal{SU}_{X}(r, \Ox) \dashrightarrow |r\Theta|$.
\par
Henceforth we suppose that $r = 2n$ and consider the map
\[ D \colon \Mn \dashrightarrow |2n\Theta|. \] 
We can be more precise about the image of $D$. Recall that the Serre duality involution $\iota \colon J^{g-1} \to J^{g-1}$ is given by $L \mapsto K_{X}L^{-1}$. Since $\iota^{*}\Theta = \Theta$, we have induced involutions $\iota^*$ on $H^{0}(J^{g-1}, 2n\Theta)$. The projectivisations of the $+1$ and $-1$ eigenspaces of this involution correspond to the spaces of $\iota^*$-invariant divisors and are denoted $|2n\Theta|_+$ and $|2n\Theta|_-$ respectively. We have $h^{0}(J^{g-1} , 2n \Theta)_{\pm} = 2n^{g} \pm 2^{g-1}$.

\begin{lemma}
The image of $D \colon \Mn \dashrightarrow |2n\Theta|$ is contained in $|2n\Theta|_+$. \label{sympdivW}
\end{lemma}
\textbf{Proof}\\
See Beauville \cite{Bea2004}, section 2.\\
\\
We denote $\Xi$ the line bundle $D^{*}O(1)$. We claim that the base locus of $|\Xi|$ is exactly the set of bundles $W$ with $S(W) = J^{g-1}$. For any $L \in J^{g-1}$, we write $H_L$ for the hyperplane of divisors containing $L$, that is, the image of $L$ under the standard map $\phi_{2n \Theta_+} \colon J^{g-1} \to |2n\Theta|_{+}^*$. Then
\begin{align*} S(W) = J^{g-1} & \Longleftrightarrow \: W \in \bigcap_{L \in J^{g-1}} \Supp ( D^{*}H_{L} ), \hbox{ by definition} \\
& \Longleftrightarrow \: W \in \bigcap_{H \in |2n\Theta|_+} \Supp ( D^{*}H ), \\
& \hspace{2cm} \hbox{since the image of $\phi_{2n\Theta_+}$ is nondegenerate} \\
& \Longleftrightarrow \: \hbox{$W$ belongs to every $\Xi$ divisor, as $|2n\Theta|_{+} = |\Xi|^*$}\\
& \Longleftrightarrow \: \hbox{$W$ is a base point of $|\Xi|$.} \end{align*}

Several results on such bundles, together with an example which we will see later, were given by Raynaud in \cite{Ray1982}. A useful survey of the area by Beauville can be found in \cite{Bea1995}. More examples of bundles without theta divisors were found by Popa in \cite{Popa1999}, and Schneider in \cite{Sch2005} gives results on the dimension of the locus of such bundles in some of the moduli spaces ${\cal U}(r, r(g-1))$.
\par
In this paper we focus on the case of semistable symplectic bundles of rank four (so $n=2$) over a curve of genus two. Here $\Mt$ is of dimension ten and $D$ is a rational map $\Mt \dashrightarrow \ftp = |\Xi|^{*}= \PP^9$. Some time ago, Arnaud Beauville found that some of Raynaud's bundles in this case admit symplectic structures, and conjectured that these were the only possibilities in this situation. The aim of this paper is to prove Beauville's conjecture:

\begin{thm} If $X$ is of genus two then the base locus of the linear system $|\Xi|$ on $\Mt$ consists of six points, which are in canonical bijection with the Weierstrass points of $X$. \label{main} \end{thm}

Here is a summary of the paper. We begin with two results for genus $g$. One describes certain isotropic subsheaves of stable symplectic bundles of rank $2n$ without theta divisors. The second is that in rank four, any base points of $|\Xi|$ must be stable vector bundles.
\par
We then specialise to $g=2$. Firstly, we show that the expected number of base points is six in this case.
\par
An important ingredient is a description of this $\Mt$ from \cite{Hit2005b} using vector bundle extensions. Let $W \to X$ be a symplectic or orthogonal bundle of rank $2n$ and $E \subset W$ an \textsl{Lagrangian} subbundle (that is, isotropic of maximal rank $n$). Then $W$ is naturally an extension of $E^*$ by $E$, and defines a class $\delta(W) \in H^{1}(X, \Homom(E^{*}, E))$. Conversely, it is natural to ask which such extensions are induced by bilinear forms. We have
\begin{thm}
Let $W$ be an extension of $E^*$ by $E$. Then $W$ carries a symplectic (resp., orthogonal) structure with respect to which $E$ is isotropic if and only if $W$ is isomorphic as a vector bundle to an extension whose class belongs to the subspace $H^{1}(X, \Sym^{2}E)$ (resp., $H^{1}(X, \bigwedge^{2}E)$ of $H^{1}(X, \Homom(E^{*}, E))$). \label{symporthext} \end{thm}
\textbf{Proof}\\
See \cite{Hit2005a}, Criterion 2.\\
\\
This motivates the following result:

\begin{thm} Let $X$ be a curve of genus two. Every semistable symplectic vector bundle of rank four over $X$ is an extension $0 \to E \to W \to E^{*} \to 0$ for some stable bundle $E \to X$ of rank two and degree $-1$, where $E$ is isotropic and the class $\delta(W)$ is symmetric. Moreover, for generic $W$, there are $24$ such $E$. Equivalently, the moduli map
\[ \Phi \colon \bigcup_{\begin{array}{c} \hbox{$E$ stable of rank two} \\ \hbox{and degree $-1$} \end{array}} \PP H^{1}(X, \Sym^{2}E) \dashrightarrow \Mt \]
is a surjective morphism\footnote{We must of course give this union of projective spaces a suitable algebraic structure.} which is generically finite of degree $24$. \label{cover} \end{thm}
\textbf{Proof}\\ This is the main result of \cite{Hit2005b}.\\
\\
We denote the five-dimensional projective space $\PP H^{1}(X, \Sym^{2}E)$ by $\pfe$. By Theorem \ref{cover}, it suffices to search for base points of $|\Xi|$ in each $\pfe$ separately. In $\S$ 5, we construct a rational map $\J \dashrightarrow \pfed$ with some useful properties, and use it to prove that no $\pfe$ contains more than one base point of $|\Xi|$.
\par
In $\S$ 6, we consider extensions $0 \to O_{X}(-w) \to E \to O_{X} \to 0$ where $w$ is a Weierstrass point of $X$. We show that this $\pfe$ contains a bundle $W_w$ without a theta divisor, and furthermore ($\S$ 7), that the isomorphism class of $W_w$ does not depend on the class of the extension $E$ in $H^{1}(X, \Homom(O_{X}, O_{X}(-w)))$.
\par
In $\S$ 8, we prove that any $W \in \Mt$ without a theta divisor must contain some such family of extensions $E$ as isotropic subbundles, so is of this form. To conclude, we show that the base locus is reduced, so consists of six points.\\
\par
In $\S$ 9, we recall a construction of Raynaud \cite{Ray1982} of stable bundles of rank four and degree zero over $X$ which have no theta divisors, and describe Beauville's construction of symplectic structures on some of these bundles. We then show how the bundles constructed in $\S$ 6 correspond to these ones.
\par
In the last section, we use Thm.\ \ref{main} to show that $D \colon \Mt \dashrightarrow \ftp$ is surjective, and notice that there exist stable bundles in $\Mt$ with reducible theta divisors.

\section{Some isotropic subsheaves}

In this section, $X$ will have genus $g \geq 2$. Firstly, we quote a couple of technical results. Let $F$ and $G$ be vector bundles over $X$. We describe two maps between associated cohomology spaces. The first one is the cup product
\[ \cup \colon H^{1}(X, \Homom(G,F)) \to \Homom( H^{0}(X, G) \to H^{1}(X, F) ), \]
and the second the natural multiplication map of sections
\[ m \colon H^{0}(X, G) \otimes H^{0}(X, K_{X} \otimes F^{*}) \to H^{0}(X, K_{X} \otimes F^{*} \otimes G). \]
\begin{prop}
The maps $\cup$ and $m$ are canonically dual, via Serre duality. \label{multdualcup}
\end{prop}
\textbf{Proof}\\
This is well known; see for example \cite{Hit2005}, chap.\ 6.\\
\\
We will also use the map $c \colon H^{1}(X, \Ox) \to H^{1}(X, \Endom(M \otimes W))$ induced by $\lambda \mapsto \lambda \cdot \Iden_{M \otimes W}$.\\
\\
\textbf{Notation:} We denote the sheaf of regular sections of a vector bundle $W$, $L$, $\Ox$, $\Kx$ by the corresponding script letter $\mathcal W$, $\mathcal L$, $\strs$, $\cans$.\\
\\
We now give a result for symplectic bundles of rank $2n$. We will need an adaptation of Prop.\ 2.6 (1) in Mukai \cite{Muk1995}. Recall that the tangent space to $J^{g-1}$ at any point is isomorphic to $H^{1}(X, \Ox)$. The number $h^{0}(X, L \otimes W)$ is constant on an open subset of $S(W)$. Let $M$ belong to this subset.

\begin{prop} The tangent space to $S(W)$ at $M$ is isomorphic to
\[ \Ker \left( \cup \circ c \colon H^{1}(X, \Ox) \to \Homom \left( H^{0}(X, M \otimes W) , H^{1}(X, M \otimes W) \right) \right). \]
\label{ZartangSW} \end{prop}
\textbf{Proof}\\
Recall that a deformation of $M \otimes W$ is an exact sequence
\[ 0 \to M \otimes W \to \mathbf{V} \to M \otimes W \to 0. \]
We are interested in those deformations which are of the form $c(v)$. By definition, $c(v)$ is tangent to the locus $S(W)$ if and only if all global sections of $M \otimes W$ lift to the extension $\mathbf{V}_{c(v)}$. But we have the cohomology sequence
\begin{multline*} 0 \to H^{0}(X, M \otimes W ) \to H^{0}( X , \mathbf{V}_{c(v)}) \to H^{0}(X, M \otimes W) \\ \to H^{1}(X, M \otimes W) \to \cdots \end{multline*}
Following Kempf \cite{Kem1983}, one shows that the boundary map is none other than cup product by $c(v)$. Therefore, \textit{all} global sections of $M \otimes W$ lift if and only if cup product by $c(v)$ is zero, that is, $v \in \Ker ( \cup \circ c )$. \qed \\
\\
Now we can prove
\begin{thm} Let $W$ be a stable symplectic bundle $W$ of rank $2n$ over $X$. Then $h^{0}(X, L \otimes W) > 0$ for every $L \in J^{g-1}$ if and only if ${\mathcal W}$ has an isotropic subsheaf ${\mathcal L}^{-1} \oplus \cans^{-1} {\mathcal L}$ for generic $L \in J^{g-1}$. \label{isotsubsh} \end{thm}
\textbf{Proof}\\
The ``if'' is clear. Conversely, suppose $h^{0}(X, L \otimes W) > 0$ for all $L \in J^{g-1}$. Let $M \in S(W)$ be such that $h^{0}(X, M \otimes W)$ is the generic value. Then the tangent space to $S(W)$ at $M$ is the whole of $H^{1}(X, \Ox)$. By Prop.\ \ref{ZartangSW}, the map $\cup \circ c$ is zero.
\par 
Now by Prop.\ \ref{multdualcup}, the map $\cup$ is dual to the multiplication
\[ H^{0}(X, M \otimes W) \otimes H^{0} ( X, \Kx M^{-1} \otimes W^{*} ) \to H^{0}(X, \Homom( \Kx \otimes W^{*} \otimes W ) ). \]
Moreover, one can show that $c$ is dual to the trace map
\[ \tr \colon H^{0}(X, \Homom(M \otimes W , \Kx \otimes ( M \otimes W ))) \to H^{0}(X, \Kx ). \]
Since $\cup \circ c$ is zero, then, so is $\tr \circ m$.
\par
Now the trace map can be identified with the map defined on decomposable elements by $e^{*} \otimes f \mapsto e^{*}(f)$. Also, $W$ is self-dual via the symplectic form $\omega$, and the induced isomorphism $W \xrightarrow{\sim} W^*$ is unique up to scalar since $W$ is stable.
\par
Combining these facts, the vanishing of $\tr \circ m$ means that the images of the maps $M^{-1} \to W$ and $\Kx^{-1}M \to W$ annihilate under contraction. Moreover, any line bundle is isotropic with respect to a symplectic form. If $M^{2} \neq \Kx$, then, there is an isotropic subsheaf ${\mathcal M}^{-1} \oplus \cans^{-1} {\mathcal M}$ in ${\mathcal W}$. This completes the proof of Thm.\ \ref{isotsubsh}. \qed \\

\section{Stability of symplectic bundles without theta divisors}

Henceforth, we assume that $n=2$. We will need
\begin{prop}
Let $X$ be a curve of genus $g \geq 2$ and $F \to X$ a semistable bundle of rank at most two and degree zero. Then $h^{0}(X, L \otimes F) = 0$ for generic $L \in J^{g-1}$. \label{rk2thetadiv} \end{prop}
\textbf{Proof}\\
This follows from Raynaud \cite{Ray1982}, Prop.\ 1.6.2.

\begin{lemma}
Let $X$ be a curve of genus $g \geq 2$. Any semistable symplectic bundles of rank four over $X$ without theta divisors are in fact \emph{stable} vector bundles. \label{bpstable} \end{lemma}
\textbf{Proof}\\
We show that every strictly semistable symplectic vector bundle $W$ over $X$ of rank four admits a theta divisor. It can be shown\footnote{This is proven in \cite{Hit2005}, chap.\ 2; the arguments are adapted from those of Ramanan \cite{Ram1981} and Ramanathan \cite{Rthn1975} for orthogonal bundles.} that such a $W$ is $S$-equivalent to a direct sum of stable bundles of rank one and/or two and degree zero. Thus it suffices to prove that every such direct sum admits a theta divisor. This follows from the last proposition. \qed

\section{The number of base points in the genus two case}

For the rest of the paper, we suppose that $X$ has genus two. In this section we find the expected number of base points of $|\Xi|$ in this case.

\subsection*{Determinant bundles}

Here we recall very briefly some facts about line bundles over the moduli space $\Mxs$. For a general treatment of this kind of question, we refer to Beauville--Laszlo--Sorger \cite{BLS1998}, Laszlo--Sorger \cite{LS1997} and Sorger \cite{Sor1998}.\\
\par
To a representation of $\Spt$, we can associate a line bundle over $\Mxs$, called the \textsl{determinant bundle} of the representation. $\Xi$ is the determinant bundle of the standard representation of $\Spt$, and the Picard group of $\Mxs$ is $\Z \cdot \Xi$. To a representation $\rho$ of $\Spt$ we associate a number $d_{\rho}$ called the \textsl{Dynkin index} of $\rho$, and the determinant bundle of $\rho$ is $\Xi^{d_{\rho}}$. The canonical bundle of $\Mxs$ is the dual of the determinant bundle of the adjoint representation, and is therefore $\Xi^{-6}$ by Sorger \cite{Sor1998}, Tableau B, since $\Spt$ is of type $C_2$.

\begin{prop}
If the base locus of $|\Xi|$ is of dimension zero then it consists of six points, counted with multiplicity.\footnote{This calculation was first done by Arnaud Beauville.} \label{count}
\end{prop}
\textbf{Proof}\\
The dimension of $\Mt$ is ten, so if the base locus of $|\Xi|$ is of dimension zero then its scheme theoretic length is given by $\mathrm{c}_{1}(\Xi)^{10}$. To calculate this number, we follow an approach of Laszlo \cite{Las1996}, $\S$ V, Lemma 5. The Hilbert function of $\Xi$ is defined as $n \mapsto \chi(\Mt,\Xi^{n})$. For large enough $n$, this coincides with a polynomial $p(n)$. We claim that the leading term of $p(n)$ is $\mathrm{c}_{1}(\Xi)^{10}/10!$. To see this, suppose $\alpha$ is the Chern root of $\Xi$. Then, by Hirzebruch--Riemann--Roch,
\[ \chi(\Mt,\Xi^{n}) = \int_{\Mt}\exp(n\alpha)\mathrm{td}(\Mt). \]
Since $\alpha^{i} = 0$ for all $i \geq 11$, the only term which contains $n^{10}$ here is $\mathrm{c}_{1}(\Xi)^{10}/10!$ as required.
\par
Now we have seen that the canonical bundle of $\Mt$ is $\Xi^{-6}$. Hence, by Serre duality, $p(n) = (-1)^{10}\chi(\Mt, \Xi^{-6-n}) = p(-6-n)$, equivalently, $p(n)$ is symmetric about $n=-3$.
\par
By a result stated on p.\ 4 of Oxbury \cite{Oxb2001a}, the spaces $H^{i}(\Mt, \Xi^{n})$ vanish for all $i > 0$ and $n > 0$. Moreover, for all $n < 0$, we have $h^{0}(\Mt, \Xi^{n}) = 0$ by stability. Thus $p(-5) =  p(-4) =  p(-3) =  p(-2) =  p(-1) = 0$. Hence $p(n)$ is equal to
\[ \gamma(n+5)(n+4)(n+3)^{2}(n+2)(n+1)(n-\alpha)(n+6+\alpha)(n-\beta)(n+6+\beta) \]
for some $\alpha, \beta, \gamma \in \mathbb{R}$. We wish to find $\gamma$. To do this, we find the values of $p$ at $0$, $1$ and $2$. By the Verlinde formula (Oxbury--Wilson \cite{OW1996}, $\S$ 2) we have $p(0) = 1$ and $p(1) = 10$ and
\[ p(2) = 2^{2} \times 5^{2} \times \sum {\mathcal S}(s,t)^{-2} \]
where 
\[ {\mathcal S}(s,t) = 2^{4} \sin\left(\frac{\pi(s+t)}{5}\right) \sin\left(\frac{\pi t}{5}\right) \sin\left(\frac{\pi s}{10}\right) \sin\left(\frac{\pi(t + 2s)}{10}\right) \]
and the sum is taken over all pairs $s,t$ with $s, t \geq 1$ and $s+t \leq 4$. 
We calculate $p(2) = 58$ with Maple. These values yield the equations
\[ \gamma \cdot 5! \times 3 (-\alpha)(6 + \alpha)(-\beta)(6 + \beta) = 1, \]
\[ \gamma \cdot 6! \times 3 (1-\alpha)(7 + \alpha)(1-\beta)(7 + \beta) = 10, \]
\[ \gamma \cdot \frac{7!}{2} \times 3 (2-\alpha)(8 + \alpha)(2-\beta)(8 + \beta) = 58. \]
Solving with Maple, we obtain $\gamma = 6 \times 10!^{-1}$, so if $\Bs|\Xi|$ is of dimension zero then it consists of six points, counted with multiplicity. \qed

\section{Study of the extension spaces}

We now begin to study the extension spaces $\pfe = \PP H^{1}(X, \Sym^{2}E)$ where $E$ is a stable bundle of rank two and degree $-1$. We describe a rational map $\J \dashrightarrow \pfed$. Firstly, we state a result similar to Thm.\ \ref{cover}:
\begin{lemma}
Every stable bundle $E \to X$ of rank two and degree $-1$ is a nontrivial extension $0 \to L^{-1} \to E \to M^{-1} \to 0$ for one, two, three or (generically) four pairs $(L,M)$ where $L$ and $M$ are line bundles of degrees one and zero respectively. Moreover, $h^{0}(X, \Homom(L^{-1}, E )) = 1 = h^{0}(X, \Homom(E, M^{-1}))$ for all such $(L,M)$. \label{lsubbsE} \end{lemma}
\textbf{Proof}\\
This follows from \cite{Hit2005b}, Lemmas 5 and 6.\\
\\
We claim that $h^{0}(X, L \otimes E^{*})\cdot h^{0}(X, K_{X}L^{-1} \otimes E^{*}) = 1$ for generic $L \in \J$. To see this, note that by Riemann--Roch, $h^{0}(X, L \otimes E^{*}) \geq 2$ if and only if $h^{1}(X, L \otimes E^{*})$ is nonzero. By Serre duality,
\[ h^{1}(X, L \otimes E^{*}) = h^{0}(X, K_{X}L^{-1} \otimes E). \]
But by Lemma \ref{lsubbsE}, this is nonzero for at most four $L$. In a similar way, we see that $h^{0}(X, K_{X}L^{-1} \otimes E^{*})$ is greater than $1$ for at most four $L$. Thus
\[ h^{0}(X, L \otimes E^{*})\cdot h^{0}(X, K_{X}L^{-1} \otimes E^{*}) \]
is different from $1$ for at most eight $L$. We write $U_E$ for the complement of these points in $\J$.
\par
For each $L \in U_E$, we can consider the composed map $\widetilde{m}$
\[ \xymatrix{ H^{0}(X, L \otimes E^{*}) \otimes H^{0}(X, K_{X}L^{-1} \otimes E^{*}) \ar[r]^-{m} & H^{0}(X, K_{X} \otimes E^{*} \otimes E^{*}) \ar[d] \\
 & H^{0}(X, K_{X} \otimes \Sym^{2}E^{*}). } \]
We claim that the image of $\widetilde{m}$ is of dimension $1$; this follows from the last paragraph and the fact that no nonzero decomposable vector is antisymmetric. Thus we can define a rational map
\[ j_{E} \colon \J \dashrightarrow \PP H^{0}(X, K_{X} \otimes \Sym^{2}E^{*}) = \pfed \]
by sending $L$ to the image of $\widetilde{m}$. This is is defined exactly on $U_E$.
\par
Now by Serre duality, a nontrivial symplectic extension $W$ of $E^*$ by $E$ with class $\delta(W) \in H^{1}(X, \Sym^{2}E)$ defines a hyperplane $H_{W} = \PP \: \Ker ( \delta(W) )$ in $\pfed$. We have
\begin{lemma}
Let $W \in \pfe$ be a symplectic extension and $L \to X$ a line bundle of degree one belonging to $U_{E} \subset \J$. Then $h^{0}(X, L \otimes W) > 0$ if and only if $j_{E}(L) \in H_W$. \label{LtensW}
\end{lemma}
\textbf{Proof}\\
Tensoring the sequence $0 \to E \to W \to E^{*} \to 0$ by $L$, we get the cohomology sequence
\[ 0 \to H^{0}(X, L \otimes W) \to H^{0}(X, L \otimes E^{*}) \xrightarrow{\cdot \cup \delta(W)} H^{1}(X, L \otimes E) \to \cdots \]
whence $h^{0}(X, L \otimes W) > 0$ if and only if the cup product map has a kernel. By hypothesis,
\[ h^{0}(X, K_{X}L^{-1} \otimes E) = h^{1}(X, L \otimes E^{*}) = 0 \]
so $h^{0}(X, L \otimes E^{*}) = 1$ by Riemann--Roch. Similarly, we see that
\[ h^{1}(X, L \otimes E) = h^{0}(X, K_{X}L^{-1} \otimes E^{*}) = 1. \]
Thus $h^{0}(X, L \otimes W) > 0$ if and only if cup product by $\delta(W)$ is zero. By Prop.\ \ref{multdualcup} (with $F = L \otimes E$ and $G = L \otimes E^*$), this means that $m^{*}\delta(W) = 0$. In other words, the image of
\[ m \colon H^{0}(X, L \otimes E^{*}) \otimes H^{0}(X, K_{X}L^{-1} \otimes E^{*}) \to H^{0}(X, K_{X} \otimes E^{*} \otimes E^{*}) \]
belongs to $\Ker(\delta(W))$. This is equivalent to $\Image ( \widetilde{m} ) \subset \Ker(\delta(W))$ because $\delta(W)$ is symmetric. Projectivising, this becomes $j_{E}(L) \in H_W$ (our hypothesis of generality on $L$ implies that $j_{E}(L)$ is defined). \qed \\

By Lemma \ref{LtensW}, we see that $\pfe$ contains an extension without a theta divisor if and only if the image of $j_{E}$ is contained in a hyperplane $H \subset \pfed$. The extension will be that $W$ such that $H_W = H$.

\begin{lemma}
For a general stable $E \to X$ of rank two and degree $-1$, the image of $j_E$ is nondegenerate in $\pfed$, and for any such $E$, it spans a hyperplane. \label{jnondeg}
\end{lemma}
\textbf{Proof}\\
By Lemma \ref{lsubbsE}, we have at least one short exact sequence
\begin{equation} 0 \to M \to E^{*} \xrightarrow{b} N \to 0 \label{E} \end{equation}
where $M$ and $N$ are line bundles of degree zero and one respectively. We claim that this induces a short exact sequence
\begin{equation} 0 \to K_{X} M^{2} \to K_{X} \otimes \Sym^{2}E^{*} \xrightarrow{c} K_{X}N \otimes E^{*} \to 0 \label{S2E} \end{equation}
where $c$ is induced by the map $E \otimes E \to N \otimes E$ given by
\[ e \otimes f \mapsto b(e) \otimes f + b(f) \otimes e .\]
We work with the map induced by $c$ on the associated locally free $\strs$-modules. Let $e, f \in {\mathcal E}_x$ be such that $e \in {\mathcal M}_x$ but $f$ is not. Then the image of $e \otimes f + f \otimes e$ belongs to ${\mathcal N \otimes M}$ but that of $f \otimes f$ does not. Thus the image contains two linearly independent elements of ${\mathcal N \otimes E}$, so the map on sheaves is surjective. The kernel of $c$ clearly contains $M^2$, so is equal to it since they are of the same rank and degree. This establishes the claim.\\
\\
Now note that the class $\langle \delta(E) \rangle \in \PP H^{1}(X, \Homom(L, M))$ can be identified with a divisor $p_{1} + p_{2} + p_{3} \in |K_{X}LM^{-1}|$. We make the following hypotheses of generality on $E$:
\begin{itemize}
\item At least one degree zero line subbundle $M \subset E^*$ is not a point of order two in $J^0$.
\item For each such $M$, there is a unique pair of points $q_{1}, q_{2} \in X$ such that
\[ K_{X}M^{2} = \strs(q_{1} + q_{2}). \]
We require that for at least one such $M \subset E$, the sets $\{ q_{1}, q_{2} \}$ and $\{ \iota p_{1}, \iota p_{2}, \iota p_{3} \}$ be disjoint.
\end{itemize}
We consider a short exact sequence (\ref{E}) where $M^2$ is nontrivial. The associated cohomology sequence is then
\[ 0 \to H^{0}(X, K_{X}M^{2}) \to H^{0}(X, K_{X} \otimes \Sym^{2}E^{*}) \to H^{0}(X, K_{X}N \otimes E^{*}) \to 0. \]
We show that $\PP H^{0}(X, K_{X}M^{2})$ is spanned by points of $j_{E}(U_{E})$. It is not hard to see that $j_{E}(L)$ is this point if and only if $L \in U_E$ and
\[ L^{-1} = M(-x) \quad \hbox{and} \quad K_{X}^{-1}L = M(-y) \]
for some points $x, y$ of the curve. This condition can be interpreted geometrically as $L \in t_{M^{-1}}\Theta \cap t_{M}\Theta$ (notice that in genus two there is a canonical isomorphism $X \xrightarrow{\sim} \Supp (\Theta)$). Since $M^2$ is nontrivial, this consists generically of $\Theta^{2} = 2$ points, which are exchanged by $\iota$. %We deduce that $K_{X}M^{2} = \strs(x+y)$; therefore
We take
\[ L = M(\iota q_{1}) = M^{-1}(q_{2}) \]
where $q_1$ and $q_2$ are as defined above; then $\Kx L^{-1} = M (\iota q_{2}) = M^{-1}(q_{1})$. We check that $j_E$ is defined at these points. It is necessary that neither $M(\iota q_{1}) = M^{-1}(q_{2})$ nor its image $M(\iota q_{2}) = M^{-1}(q_{1})$ under the Serre involution be a quotient of $E^*$ or the Serre image of a quotient of $E^*$. It is not hard to show\footnote{See for example the proof of \cite{Hit2005b}, Lemma 6.} that the degree one line bundle quotients of $E^*$ are $L^{-1}$, $M^{-1}(p_{1})$, $M^{-1}(p_{2})$ and $M^{-1}(p_{3})$. Thus we must check that
\[ \{ M(\iota q_{1}), M(\iota q_{2}) \} \cap \{ N, KN^{-1}, M(p_{i}), K_{X}M^{-1}(-p_{i}) \} = \emptyset. \]
Since the set where $j_E$ is not defined is $\iota$-invariant, it suffices to check that neither $M(\iota q_{1})$ nor $M(\iota q_{2})$ is equal to $N$ or $M(p_{i})$.
\par
If $M(\iota q_{j}) = N$ then $\PP H^{1}(X, \Homom(N, M)) \cong |K_{X}(\iota q_{j})|$ and
\[ p_{1} + p_{2} + p_{3} = \iota q_{j} + r + \iota r \]
for some $r \in X$. Thus $q_{j} = \iota p_i$ for some $i, j$. On the other hand, if $M(\iota q_{j}) = M(p_{i})$ then $q_{j} = \iota p_i$. By our second assumption of generality, then, $j_{E}$ is defined at both $M(\iota q_{j})$, and $\PP H^{0}(X, K_{X}M^{2})$ belongs to $j_{E}(U_{E})$.
\par
Now $K_X$ tensored with (\ref{E}) yields the cohomology sequence
\[ 0 \to H^{0}(X, K_{X}M) \to H^{0}(X, K_{X} \otimes E^{*}) \to H^{0}(X, K_{X}N) \to 0 \]
and we have the diagram
\[ \xymatrix{ & 0 \ar[d] & \\
 & H^{0}(X, K_{X}NM) \ar[d] & \\
 H^{0}(X, K_{X} \otimes \Sym^{2}E^{*})  \ar[r]^{c} \ar[dr]^{d} & H^{0}(X, K_{X}N \otimes E^{*})  \ar[r] \ar[d] & 0 \\
 & H^{0}(X, K_{X}N^{2}) \ar[d] & \\
 & 0 & } \]
We must show that the image of $\PP c \circ j_{E}$ is nondegenerate in $\PP H^{0}(X, K_{X}N \otimes E^{*})$. For this, it suffices to show that $\PP H^{0}(X, K_{X}NM)$ is spanned and that $\PP d \circ j_{E}(U_{E})$ is nondegenerate in $\PP H^{0}(X, K_{X}N^{2})$.
\par
By the definition of $c$, we have $\PP c \circ j_{E}(L) \in \PP H^{0}(X, K_{X}NM)$ if $L^{-1} = M(-x)$ for some $x \in X$ but $K_{X}^{-1}L \neq M(-y)$ for any $y \in X$, or vice versa. This is equivalent to $L$ belonging to the symmetric difference of $t_{M^{-1}}\Theta$ and $t_{M}\Theta$. Since $M \neq M^{-1}$, we can find infinitely many such $L$. We need to find two which define different divisors in $|K_{X}NM|$. We observe that if $x$ is not a base point of $|K_{X}NM|$ then $\PP c \circ j_{E}(M^{-1}(x))$ is the divisor in $|K_{X}NM|$ containing $x$. Thus if neither $x$ nor $y$ is a base point and they belong to different divisors of $|K_{X}NM|$ then the images of $j_{E}(M^{-1}(x))$ and $j_{E}(M^{-1}(y))$ generate $\PP H^{0}(X, K_{X}NM)$.
\par
Next, we have $\PP d \circ j_{E}(L) \in \PP H^{0}(X, K_{X}N^{2})$ if and only if neither $L^{-1}$ nor $K_{X}^{-1}L$ is of the form $M(-x)$, equivalently, $L$ does not belong to $t_{M^{-1}}\Theta \cup t_{M}\Theta$. In fact $\PP d \circ j_{E}$ is dominant. Let $x_{1} + x_{2} + x_{3} + x_{4}$ be any divisor in $|K_{X}N^{2}|$. For generic $x_1$ and $x_2$, we know that $N^{-1}(x_{1}+x_{2})$ will not belong to $t_{M^{-1}}\Theta \cup t_{M}\Theta$, so we can put $L = N^{-1}(x_{1}+x_{2})$ and then
\[ K_{X}L^{-1} = K_{X}N(-x_{1}-x_{2}) = N^{-1}(x_{3}+x_{4}). \]
Thus $\PP d \circ j_{E}$ is dominant (and generically of degree six). In particular, the image spans $|K_{X}N^{2}| = \PP^2$.\\
\par
Finally, we show that the image of $j_E$ always spans a $\PP^4$. If $E^*$ fits into a short exact sequence $0 \to M \to E^{*} \to N \to 0$ where $M$ is of order two, then we form as before a cohomology diagram
\[ \xymatrix{ & 0 \ar[d] & \\
 & H^{0}(X, K_{X}NM) \ar[d] & \\
 H^{0}(X, K_{X} \otimes \Sym^{2}E)  \ar[r]^{c} \ar[dr]^{d} & H^{0}(X, K_{X}N \otimes E)  \ar[r] \ar[d] & \cdots \\
 & H^{0}(X, K_{X}N^{2}) \ar[d] & \\
 & 0 & } \]
Since $M = M^{-1}$, we have $t_{M^{-1}}\Theta = t_{M}\Theta$ and $j_E$ is defined at infinitely many points of $t_{M^{-1}}\Theta \cap t_{M}\Theta$. Choose two distinct points $p,q \in X$ such that $p + q$ is not a canonical divisor; then $j_{E}(M(p))$ and $j_{E}(M(q))$ span $\PP H^{0}(X, K_{X}M^{2})$. Since the symmetric difference of $t_{M^{-1}}\Theta$ and $t_{M}\Theta$ is empty, there are no points of the image mapping to $\PP H^{0}(X, K_{X}NM) \subset \PP H^{0}(X, K_{X}N \otimes E)$, However, the composed map to $\PP H^{0}(X, K_{X}N^{2})$ is dominant, by the same argument. Thus we can find five linearly independent points of the image of $j_E$, spanning a $\PP^4$ in $\pfed$.
\par
On the other hand, if the second hypothesis of generality fails for some exact sequence $0 \to M \to E^{*} \to N \to 0$ where $M$ is not of order two, then $j_E$ is not defined at either $M(\iota q_{j})$. However, the rest of the proof goes through and again we can find five independent points of the image of $j_{E}(U_{E})$. \qed \\
\\
\textbf{Remark:} In fact the second generality hypothesis can be weakened slightly (we then have to blow up $\J$ at a point) but the statement as given is strong enough for our purposes.  

\section{An example of a base point}

In this section we give an explicit example of a base point of $|\Xi|$. We begin by constructing an $E$ which violates both of the generality conditions stated in the last lemma.
\par
Let $w \in X$ be a Weierstrass point. Extensions
\begin{equation} 0 \to O_{X}(-w) \to E \to O_{X} \to 0 \label{specialE} \end{equation}
are determined as vector bundles by a divisor in
\[ |K_{X}(w)| \cong \PP H^{1}(X, \Homom(O_{X}, O_{X}(-w))). \]
This divisor will be of the form $w + p + \iota p$ where $\iota$ is the hyperelliptic involution on $X$. It is not hard to see that the degree $-1$ subbundles of $E$ are $O_{X}(-w)$, $O_{X}(-p)$ and $O_{X}(- \iota p)$.
\begin{prop}
The bundle $E$ is $\iota$-invariant.
\end{prop}
\textbf{Proof}\\
By pulling back the sequence (\ref{specialE}) by $\iota$, we get a short exact sequence
\[ 0 \to \Ox (-w) \to \iota^{*} E \to \Ox \to 0 \]
since $\Ox(-w)$ and $\Ox$ are $\iota$-invariant. The divisor determining the extension $\iota^{*}E$ is $\iota(w) + \iota(p) + \iota(\iota p) = w + \iota p + p$, so $\iota^{*}E$ is isomorphic to $E$. \qed \\

We now give some more information on $D$. For any stable $E \to X$ of rank two and degree $-1$, we have a map $\De \colon \pfe \dashrightarrow \ftp$ which is the composition of $D \colon \Mt \dashrightarrow \ftp$  with the classifying map $\Phi|_{\pfe} \colon \pfe \to \Mt$.
\begin{lemma} The map $\De \colon \pfe \dashrightarrow \ftp$ is linear. \label{DElin} \end{lemma}
\textbf{Proof}\\
The degree of $\De$ is constant with respect to $E$ since the moduli space of such $E$ is connected (see for example Le Potier \cite{LeP1997}). It therefore suffices to prove the lemma for a general $E$. In particular, we can suppose that $\De$ is defined everywhere on $\pfe$.
\par
For $L \in \J$, let $H_{L} \subset \ftp$ denote the hyperplane of divisors containing $L$. By definition, $D_{E}\langle \delta(W) \rangle$ belongs to $H_L$ if and only if $h^{0}(X, L \otimes W) > 0$. We show that this is a linear condition on $\pfe$. Let $W$ be an extension of $E^*$ by $E$. We get a cohomology sequence
\begin{multline*}
0 \to H^{0}(X, L \otimes W) \to H^{0}(X, L \otimes E^{*}) \xrightarrow{\cup_{L} \delta(W)} \\ H^{1}(X, L \otimes E) \to H^{1}(X,L \otimes W) \to 0.
\end{multline*}
If $h^{0}(X, L \otimes E^{*}) \geq 2$ then $h^{0}(X, K_{X} L^{-1} \otimes E) \geq 1$ and the image of $\De$ is contained in $H_{K_{X} L^{-1}} = H_L$, so $\De^{*}H_L$ is the whole space.
\par
If $h^{0}(X, L \otimes E^{*}) = 1$ then $H^{0}(X, L \otimes W)$ is nonzero if and only if $\cup_{L} \delta(W)$ is the zero map. Thus, set-theoretically, $\De^{*}H_{L} = \cup_{L}^{-1}(0)$. Now $\cup_L$ is a linear map
\[ H^{1}(X, \Sym^{2}E^{*}) \to \Homom( H^{0}(X, L \otimes E), H^{1}(X, L \otimes E)). \]
Since the latter space has dimension 1, $\Ker(\cup_{L})$ is a hyperplane.\\
\\
It remains to show that the multiplicity of $\De^{*}H_L$ is 1. We begin by showing that $\De(\pfe)$ is contained in a $\PP^5$ in $\ftp$. By generality of $E$, we can suppose that $E$ has four distinct line subbundles of degree $-1$ not including any pairs of the form $N, \Kx^{-1} N^{-1}$. Thus the divisors in the image of $\De$ all contain four points distinct modulo $\iota$. We show that these impose independent conditions on $|4\Theta|_+$.
\par
Let $L^{-1} \subset E$ be a line subbundle of degree $-1$ and $M^{-1}$ the quotient of $E$ by $L^{-1}$. Then $E$ is defined by a divisor $p_{1} + p_{2} +p_{3} \in |\Kx M L|$. By generality, we can suppose that this linear system is base point free and that the $p_i$ are distinct. The other degree $-1$ subbundles of $E$ are then $M^{-1}(-p_{1})$, $M^{-1}(-p_{2})$ and $M^{-1}(-p_{3})$. We construct even $4\Theta$ divisors $D_{0}, D_{1}, \ldots, D_4$ containing none, one, two, three and all of these points respectively. Since $\phi_{4\Theta_+} \colon \J \to \ftp^*$ descends to an embedding of the Kummer variety and the points $L, M(p_{i})$ are distinct on the Kummer, we can easily find $D_0$ and $D_1$. For $D_2$, we use the fact that $\phi_{2\Theta}$ also gives an embedding of the Kummer. Thus we can find a $2\Theta$ divisor $G$ containing $M(p_{1})$ but none of other points, and another, $G^{\prime}$, containing $M(p_{2})$ but none of the others. Since every $2\Theta$ divisor is even, the sum $G + G^{\prime}$ is an even $4\Theta$ divisor containing exactly two of the points. For $D_3$, we take $2(t_{M} \Theta + t_{M^{-1}} \Theta)$. By generality, $L$ is not of the form either $M(x)$ or $M^{-1}(x)$ for any $x \in X$, so this divisor contains the three $M(p_{i})$ but not $L$. Finally, choose any symplectic extension $W$ of $E^*$ by $E$ which has a theta divisor $D(W)$; these exist by Lemma \ref{jnondeg}. Then we take $D_{4} = D(W)$; this contains all the points.
\par
The four degree $-1$ line subbundles of $E$ thus impose independent conditions on the divisors in $\ftp = \PP^9$. Hence the image of $D_E$ is contained in a $\PP^5$. By generality and Lemma \ref{jnondeg}, the map $\De$ is a morphism $\PP^{5} \to \PP^5$, so must be surjective, therefore a finite cover. If not an isomorphism, it must be branched over a hypersurface. Now the image of $\phi_{4\Theta_+} \colon \J \to \ftp^*$ is nondegenerate, so we can find an $H_L$ whose intersection with the image of $\De$ is not contained in this branch locus. Then $\De^{*}H_L$ is reduced, and we have seen that its support is a hyperplane.\\
\\
This completes the proof of Lemma \ref{DElin}.   \qed \\

Now since $E$ is $\iota$-invariant, we can lift $\iota$ to linearisations $\tilde{\iota}$ and $\Sym^{2}(\tilde{\iota})$ of $E$ and $\Sym^{2}E$. The latter acts on $H^{1}(X, \Sym^{2}E)$, taking the class of an extension $W$ to that of $\iota^{*}W$ (modulo a scalar). Since this is an involution, we can decompose $H^{1}(X, \Sym^{2}E)$ into $+1$ and $-1$ eigenspaces $H^{1}(X, \Sym^{2}E)_{\pm}$.\\
Suppose now that $E$ is of the form (\ref{specialE}). Let
\[ \widetilde{\De} \colon H^{1}(X, \Sym^{2}E) \to H^{0}( \J , 4 \Theta )_+ \]
be a linear lift of $\De$. The kernel of $\widetilde{\De}$ consists of exactly the extensions without theta divisors. Now when $D(W)$ exists, we have $\iota^{*}D(W) = D(\iota^{*}W)$ since their supports are equal and they are both even $4 \Theta$ divisors. Therefore, for all $W \in \pfe$, either
\[ \widetilde{\De}( \delta(\iota^{*}W) ) = \widetilde{\De} ( \delta(W) ) \quad \hbox{or} \quad \widetilde{\De}( \delta(\iota^{*}W) ) = - \widetilde{\De} (\delta(W) ). \]
This means that one of the spaces $H^{1}(X, \Sym^{2}E)_{\pm}$ belongs to the kernel of $\widetilde{D_E}$. We calculate the dimensions of these eigenspaces. The key tool is the fixed point formula of Atiyah and Bott:
\begin{thm}
Let $M$ be a compact complex manifold and $\gamma \colon M \to M$ an automorphism with a finite set $\Fix(\gamma)$ of fixed points. Suppose that $\gamma$ lifts to a linearisation $\widetilde{\gamma}$ of a holomorphic vector bundle $V \to M$. Then
\[ \sum_{j} (-1)^{j} \tr \left( \gamma |_{H^{j}(M, V)} \right) = \sum_{p \in \Fix ( \gamma )} \frac{\tr \left( \gamma|_{E|_p} \right)}{\det \left( I - d\gamma_{p} \right) }. \] \label{ABLfpf} \end{thm}
\textbf{Proof}\\
See Atiyah--Bott \cite{AB1968}, Theorem 4.12.\\
\\
For us $M = X$, with $\gamma = \iota$, and $V = \Sym^{2}E$. We write down the linearisation $\tilde{\iota}$ of $E$ explicitly and determine the action of $\Sym^{2}(\tilde{\iota})$ on the fibre of $\Sym^{2}E$ over the Weierstrass points. For any bundle $V$ with sheaf of sections ${\mathcal V}$, the fibre of $V$ at a point $x$ is identified with ${\mathcal V}_{x}/m_{x}{\mathcal V}_x$ where $m_x$ is the maximal ideal of the ring ${\mathcal O}_{X,x}$ (see for example Le Potier \cite{LeP1997}, chap.\ 1). Let $z$ be a uniformiser at $w$. Then $\iota$ is given near $w$ by $z \mapsto -z$.
\par
Now it is not hard to show that a linearisation of an involution on a line bundle is unique up to multiplication by $-1$. Modulo this, the last paragraph shows:
\begin{itemize}
\item A lift of $\iota$ to $\Ox$ acts trivially on $\Ox|_p$ for all Weierstrass points $p$.
\item Since $\Ox(-w)|_{p} \cong \Ox|_p$ for each $p \neq w$, the same is true for a lift of $\iota$ to $\Ox(-w)$ at $\Ox(-w)|_p$ for these $p$.
\item A lift of $\iota$ to $\Ox(-w)$ acts by $-1$ on $\Ox(-w)|_w$.
\end{itemize}
\par
Near the point $w$, the bundle $E$ looks like $\Ox(-w) \oplus \Ox$. This is far from canonical, but we are interested only in the trace of a linearisation, which is independent of the trivialisation. We normalise $\tilde{\iota}$ such that the induced linearisation on $\Ox(-w)$ acts on the fibre over a Weierstrass point $p$ by
\[ \left\{ \begin{array}{l} -1 \hbox{ if } p \neq w \\ +1 \hbox{ if } p=w. \end{array} \right. \]
Then $\tilde{\iota}$ acts on the fibres $E|_p$ by either
\[ \left\{ \begin{array}{l} \begin{pmatrix} -1 & 0 \\ 0 & -1 \end{pmatrix} \hbox{ if } p \neq w \\ \begin{pmatrix} 1 & 0 \\ 0 & -1 \end{pmatrix} \hbox{ if } p = w \end{array} \right. \quad \hbox{or} \quad \left\{ \begin{array}{l} \begin{pmatrix} -1 & 0 \\ 0 & 1 \end{pmatrix} \hbox{ if } p \neq w \\ \begin{pmatrix} 1 & 0 \\ 0 & 1 \end{pmatrix} \hbox{ if } p = w. \end{array} \right. \]
Using the sequence $0 \to E(-w) \to \Sym^{2}E \to \Ox \to 0$, which is derived from (\ref{S2E}), we find that $\Sym^{2}(\tilde{\iota})$ acts as follows on $\Sym^{2}E|_p$:
\begin{equation} \left\{ \begin{array}{l} \diag(1,1,1) \hbox{ if } p \neq w \\ \diag(1,-1,1) \hbox{ if } p = w \end{array} \right. \quad \hbox{or} \quad \left\{ \begin{array}{l} \diag(1,-1,1) \hbox{ if } p \neq w \\ \diag(1,1,1) \hbox{ if } p = w. \end{array} \right. \label{2poss} \end{equation}
Suppose the first possibility occurs. Since $\iota$ acts by $z \mapsto -z$ in a neighbourhood of a Weierstrass point, $d\iota|_{p} = -1$. Thm.\ \ref{ABLfpf} then gives
\[ h^{1}(X, \Sym^{2}E)_{-} - h^{1}(X, \Sym^{2}E)_{+} = \frac{ 5 \times 3 + 1 \times 1}{2} = 8 \]
which is impossible since $h^{1}(X, \Sym^{2}E)_{-} + h^{1}(X, \Sym^{2}E)_{+} = 6$.
\par
Therefore, the second possibility in (\ref{2poss}) occurs, and Thm.\ \ref{ABLfpf} gives
\[ h^{1}(X, \Sym^{2}E)_{-} - h^{1}(X, \Sym^{2}E)_{+} = \frac{ 5 \times 1 + 1 \times 3}{2} = 4. \]
Solving, we obtain $h^{1}(X, \Sym^{2}E)_{+} = 1$ and $h^{1}(X, \Sym^{2}E)_{-} = 5$.
\par
By Lemma \ref{jnondeg}, the kernel of $\widetilde{\De}$ is of dimension at most one, so it must be equal to $H^{1}(X, \Sym^{2}E)_+$.\\
\\
\textit{Thus, the point $\PP H^{1}(X, \Sym^{2}E)_{+} \in \pfe$ defines a stable rank four symplectic bundle without a theta divisor, that is, a base point of $|\Xi|$ in $\Mt$.}\\
\par
In the following sections, we will show that in fact every base point of $|\Xi|$ is of this form.

\section{Maximal Lagrangian subbundles}

We write $E_e$ for the extension of $\Ox$ by $O_{X}(-w)$ defined by
\[ e \in \PP H^{1}(X, \Homom(\Ox , \Ox(-w))) = \PP^{1}. \]
We denote $W$ the base point of $|\Xi|$ constructed from $E_e$ in the last section.

\begin{lemma} The isomorphism class of the bundle $W$ is independent of the extension class $e$. \label{Windepe} \end{lemma}
\textbf{Proof}\\
We show firstly that $W$ contains a subbundle of the form $E_f$ for \textit{every} $f \in \PP H^{1}(X, \Homom(\Ox , \Ox(-w)))$.
\par
Since we have homomorphisms $\beta \colon E_{f} \to \Ox$ and $\gamma \colon \Ox \to E_{e}^*$, we can find a map $\gamma \circ \beta \colon E_{f} \to E_{e}^*$ for every $f$. Since $h^{1}(X, \Homom( E_{f}, \Ox)) = 1$, the map $\beta$ must be equivariant or antiequivariant; examining the action of the linearisation induced on $\Homom( E_{f}, \Ox)$ by those on $E_f$ and $\Ox$, we find that it is invariant. We check similarly that $\gamma$ is equivariant.
\par
We show that $\gamma \circ \beta$ factorises via $W$. By Narasimhan--Ramanan \cite{NR1969}, Lemma 3.3, this happens if and only if $\delta(W)$ lies in the kernel of the map $(\gamma \circ \beta)^{*} \colon H^{1}(X, \Sym^{2}E_{e}) \to H^{1}(X, \Homom(E_{f} , E_{e}))$. This map factorises as
\[ H^{1}(X, \Sym^{2}E_{e}) \xrightarrow{\gamma^*} H^{1}(X, \Homom(\Ox , E_{e})) \xrightarrow{\beta^*}  H^{1}(X, \Homom(E_{f} , E_{e})). \]
We see that $\gamma^{*} (\delta(W) )$ is not zero because then, by the same result, $\Ox$ would be a subbundle of $W$, which is excluded by Lemma \ref{bpstable}. Now the space $H^{1}(X, \Homom(\Ox , E_{e}))$ also has an action of $\iota$, for which $\gamma^*$ is equivariant. Thus it will be enough to show that $\beta^{*}(H^{1}(X, \Homom(\Ox , E_{e}))_{+})$ is zero.
\par
Taking $\Homom( - , E_{e})$ of the sequence $O \to \Ox(-w) \to E_{f} \xrightarrow{\beta} \Ox \to 0$, we find that $\beta^*$ fits into the exact sequence
\begin{multline*} \cdots \to H^{1}(X, \Homom(\Ox , E_{e})) \xrightarrow{\beta^*} H^{1}(X, \Homom(E_{f} , E_{e})) \to \\ H^{1}(X, \Homom(\Ox(-w) , E_{e})) \to 0. \end{multline*}
Now all the maps on cohomology are $\iota$-equivariant, so this sequence splits into a direct sum of invariant and antiinvariant sequences. We calculate the numbers $h^{1}(X, \Homom(E_{f} , E_{e}))_+$ and $h^{1}(X, \Homom(\Ox(-w) , E_{e}))_+$.
We can assume $e \neq f$, so $h^{0}(X, \Homom(E_{f}, E_{e})) = 0$ by stability, and so
\[ h^{1}(X, \Homom(E_{f}, E_{e}))_{+} + h^{1}(X, \Homom(E_{f}, E_{e}))_{-} = 4 \]
by Riemann--Roch. Locally, $\Homom(E_{f} , E_{e})$ splits into
\begin{multline*} \Homom( \Ox(-w), \Ox(-w) ) \oplus \Homom( \Ox, \Ox(-w) ) \\ \oplus \Homom( \Ox(-w), \Ox ) \oplus \Homom( \Ox, \Ox ) \end{multline*}
so, with linearisation as on $E$ in the last section, the action of a lifting of $\iota$ to $\Homom(E_{f}, E_{e})$ on the fibre over a Weierstrass point $p$ is given by
\[ \left\{ \begin{array}{ll} \diag(1,-1,-1,1) & \hbox{ if } p \neq w \\ \diag(1,1,1,1) & \hbox{ if } p = w. \end{array} \right. \]
Then Theorem \ref{ABLfpf} gives us $h^{1}(X, \Homom(E_{f}, E_{e}))_{-} - h^{1}(X, \Homom(E_{f}, E_{e}))_{+} = 2$, so $h^{1}(X, \Homom(E_{f}, E_{e}))_{-} = 3$ and $h^{1}(X, \Homom(E_{f}, E_{e}))_{+} = 1$.
\par
As for $h^{1}(X, \Homom(\Ox(-w) , E_{e}))_+$: the bundle $\Homom( \Ox(-w), E_{e})$ has global sections, of which we must take account when using Theorem \ref{ABLfpf}. There is a single independent section $\Ox(-w) \to E_e$ by the second statement of Lemma \ref{lsubbsE}, which, as before, is equivariant. Taking this into account, we calculate, using Theorem \ref{ABLfpf} as before, that
\[ h^{1}(X, \Homom( \Ox(-w), E_{e}))_{-} - h^{1}(X, \Homom( \Ox(-w), E_{e}))_{+} = 0. \]
Then by Riemann--Roch we have $h^{1}(X, \Homom( \Ox(-w), E_{e})) = 2$, and hence $h^{1}(X, \Homom( \Ox(-w), E_{e}))_{+} = 1$.
\par
Putting all this together, we have an exact sequence of vector spaces
\begin{multline*} H^{1}(X, \Homom( \Ox, E_{e} ))_{+} \xrightarrow{\beta^*} H^{1}(X, \Homom(E_{e} , E_{f} ))_{+} \\ \to H^{1}(X, \Homom( \Ox(-w), E_{e})_{+} \to 0. \end{multline*}
where each of the last two spaces is of dimension one. Thus, in particular $\beta^{*}(H^{1}(X, \Homom( \Ox, E_{e} ))_{+}) = 0$, so $E_f$ belongs to $W$, as we wanted.\\
\par
Now we have a diagram
\[ \xymatrix{ & 0 \ar[d] & 0 \ar[d] & 0 \ar[d] & \\
0 \ar[r] & \Ox(-w) \ar[d]^{\gamma} \ar[r] & E_{f} \ar[d] \ar[r] & \Ox \ar[r] \ar[d] & 0 \\
0 \ar[r] & E_{e} \ar[d] \ar[r]^{\beta} & W \ar[d] \ar[r]^{^{t}\beta \circ \omega^{-1}} & E_{e}^{*} \ar[r] \ar[d]^{^{t}\gamma} & 0 \\
0 \ar[r] & \Ox \ar[d] \ar[r] & F \ar[d] \ar[r] & \Ox(w) \ar[r] \ar[d] & 0 \\
 & 0 & 0 & 0 & } \]
where $\omega \colon W \to W^*$ is an isomorphism, unique up to scalar. The only maps whose existences are not immediate are those $\Ox \to F$ and $F \to \Ox(w)$. For the first: since $E_{e} \cap E_f$ is the subbundle $\Ox(-w)$ of $W$, the image of $E_e$ in $F$ is $E_{e}/\Ox(-w)$, which is $\Ox$. For the second, we note that $\det(F) = \det(W) \det(E_{f})^{-1} = \Ox (w)$, so $F/\Ox$ is isomorphic to $\Ox(w)$. By the second statement of Lemma \ref{lsubbsE}, there is only one independent map each from $F \to \Ox(w)$ and $E_{e}^{*} \to \Ox(w)$, so the bottom square is commutative after multiplication of the map $F \to \Ox(w)$ by a scalar.

\begin{prop}
The bundle $F$ is isomorphic to $E_{f}^*$.
\end{prop}
\textbf{Proof}\\
Firstly, note that $F$ is not a split extension since then we would have a nonzero map $W \to \Ox$, contradicting stability.
\par
Since the map $\beta \circ \gamma \colon \Ox(-w) \to W$ factorises via $E_f$, the class $\delta(E_{f})$ belongs to the kernel of the induced map
\[ \left( \beta \circ \gamma \right)_{*} \colon H^{1}(X, \Homom(\Ox , \Ox(-w))) \to H^{1}(X, \Homom( \Ox, W)). \]
Similarly, since $\left( {^{t}\gamma} \circ {^{t}\beta} \circ \omega^{-1} \right) \colon W \to \Ox(w)$ factorises via $F$, the class $\delta(F)$ belongs to the kernel of the induced map
\[ \left( {^{t}\gamma} \circ {^{t}\beta} \circ \omega^{-1} \right)^{*} \colon H^{1}(X, \Homom( \Ox(w), \Ox)) \to H^{1}(X, \Homom(W, \Ox)). \]
Now there is a commutative diagram
\[ \xymatrix { H^{1}(X, \Homom( \Ox, \Ox(-w))) \ar[r]^{\sim} \ar[d]_{\gamma_*} & H^{1}(X, \Homom( \Ox(w), \Ox)) \ar[d]^{^{t}\gamma^*} \\
H^{1}(X, \Homom( \Ox , E_{e})) \ar[r]^{\sim} \ar[d]_{\beta_*} & H^{1}(X, \Homom( E_{e}^{*}, \Ox )) \ar[d]^{(^{t}\beta \circ \omega^{-1})^*} \\
H^{1}(X, \Homom( \Ox , W )) \ar[r]^{\sim} & H^{1}(X, \Homom( W, \Ox )) } \]
where the horizontal arrows are induced by the transposes (the lowest one factorises as
\begin{multline*} H^{1}(X, \Homom( \Ox , W)) \xrightarrow{\mathrm{transpose}} H^{1}(X, \Homom( W^{*}, \Ox )) \\
 \xrightarrow{(\omega^{-1})^*} H^{1}(X, \Homom( W , \Ox )) \quad ). \end{multline*}
We have seen that $\delta(E_{f})$ and $\delta(F)$ belong to the kernels of the composed vertical maps. Since these kernels are each of dimension 1, we see that $\delta(F)$ is proportional to $^{t}(-\delta(E_{f})) = \delta(E_{f}^{*})$, whence the proposition. \qed \\
\par
The last step is to show that the class $\delta_{f}(W) \in H^{1}(X, \Homom(E_{f}^{*}, E_{f}))$ is symmetric. We recall that the involution $\delta \mapsto -^{t}\delta$ on the extension space $H^{1}(X, \Homom(E_{f}^{*}, E_{f}))$ sends the class of an extension $V$ to that of $V^*$. 
\par
Now by Lemma \ref{Fpdesc} from the Appendix (see also the proof of Lemma \ref{twentythree}), we have $h^{0}(X, W(x)) = 1$ for all $x \in X$. A generic $E_f$ has two line subbundles of the form $\Ox(-x)$ and $\Ox(- \iota x)$ for some $x \in X$, and these generate $E_f$. Therefore we have $h^{0}(X, \Homom(E_{f}, W)) = 1$ for such an $E_f$. Since $W \cong W^*$, we deduce that there is an isomorphism of exact sequences
\[ \xymatrix{ 0 \ar[r] & E_{f} \ar[r] \ar[d]^{\wr} & W \ar[r] \ar[d]^{\wr} & E_{f}^{*} \ar[r] \ar[d]^{\wr} & 0 \\
 0 \ar[r] & E_{f} \ar[r] & W^{*} \ar[r] & E_{f}^{*} \ar[r] & 0 } \]
with classes $\delta_{f}(W)$ and $-^{t}\delta_{f}(W)$. Since the bundles $E$ and $E^*$ are simple, these two classes are proportional, so $\delta_{f}(W)$ belongs to either $H^{1}(X, \Sym^{2}E_{f})$ or $H^{1}(X, \bigwedge^{2}E_{f})$. If it were the latter then $W$ would have an orthogonal structure by Theorem \ref{symporthext}, which would contradict the stability of $W$. Thus $\delta_{f}(W) \in H^{1}(X, \Sym^{2}E_{f})$ and $W$ is the base point of $|\Xi|$ associated to $E_f$, which is unique by Lemma \ref{jnondeg}. Since the symmetry of $\delta_{f}(W)$ is a closed condition on $f$, we see that $\delta_{f}(W)$ is symmetric for all $f$.\\
\\
This completes the proof of Lemma \ref{Windepe}. \qed \\
\quad \\
In this way we associate to each Weierstrass point $w \in X$ a base point of $|\Xi|$. We denote this bundle $W_w$.

\section{Characterisation of the base points}

In this section we will prove Theorem \ref{main}. We will show that an arbitrary $W \in \Mt$ with no theta divisor must be of the form $W_w$ for some Weierstrass point $w \in X$, and that for distinct $w, v \in X$, the bundles $W_w$ and $W_v$ are mutually nonisomorphic.\\
\par
Let $W \in \Mt$, then, be a base point of $|\Xi|$. Since $S(W) = \J$, in particular we have $h^{0}(X, W(x) ) > 0$ for all $x \in X$. Now by Riemann--Roch, we have $\chi(X, W \otimes \Kx) = 4$. By Serre duality $h^{1}(X, W \otimes \Kx ) = h^{0}(X, W^{*})$, which is zero since $W$ is stable by Prop.\ \ref{bpstable}. Thus $h^{0}(X, W \otimes \Kx ) = 4$. On the other hand,
\[ h^{0}(X, W \otimes \Kx(-x)) = h^{0}(X, \Homom(\Ox(-\iota x), W)) > 0 \]
for all $x \in X$, by hypothesis. This means that the rank 4 bundle $W \otimes \Kx$ is not generated by its global sections. In other words, the evaluation map
\[ \ev \colon \Ox \otimes H^{0}(X, W \otimes \Kx) \to W \otimes \Kx \]
is not of maximal rank. We denote ${\mathcal F}$ the subsheaf of ${\mathcal W}$ corresponding to the image of $\ev$.
\begin{lemma} The subsheaf ${\mathcal F}$ of ${\mathcal W} \otimes \cans$ corresponds to a vector subbundle $F \subset W \otimes \Kx$ which has rank three and degree five, and is stable. \label{Fpdesc} \end{lemma}
\textbf{Proof}\\
This proof is straightforward but rather long. We relegate it to the appendix, in order not to interrupt the story.\\
\par
By for example taking global sections of $0 \to F \to W \otimes \Kx \to \cdots,$ we have $H^{0}(X, F) \cong H^{0}(X, W \otimes \Kx)$. Denote $L^{-1} = \Ker ( \ev )$. There is an exact sequence
\begin{equation} 0 \to L^{-1} \to \Ox \otimes H^{0}(X, F) \to F \to 0. \label{defFL} \end{equation}
Henceforth we denote $F$ by $F_L$ to emphasise that $F$ depends on $L$. Note that $\det ( F_{L} ) = L$.\\
\\
Twisting by $\Kx^{-1}$, we have shown that every $W \in \Mt$ with no theta divisor contains a rank three subbundle of the form $F_{L} \otimes \Kx^{-1}$, for some $L \in J^{5}_X$. We will now say more about the structure of $F_L$.

\begin{prop}
There is an exact sequence
\[ 0 \to \Kx^{-2} L \to F_{L} \to \Kx \oplus \Kx \to 0. \]
\label{defFLt} \end{prop}
\textbf{Proof}\\
Let us show firstly that $h^{0}(X, \Homom( F_{L}, \Kx) ) = 2$. We have an exact sequence
\[ 0 \to H^{0}(X, F_{L}^{*} \otimes \Kx) \to  H^{0}(X, \Kx) \otimes H^{0}(X, L) \xrightarrow{\mu} H^{0}(X, \Kx L ) \to \cdots \]
given by tensoring the dual of (\ref{defFL}) by $\Kx$ and taking global sections. Thus $H^{0}(X, F_{L}^{*} \otimes \Kx ) = \Ker (\mu)$. Now since $X$ is of genus two, it is not hard to see that there is a short exact sequence
\[ 0 \to \Kx^{-1} \to \Ox \otimes H^{0}(X, \Kx ) \to \Kx \to 0. \]
Tensoring by $L$ and taking cohomology, we obtain
\[ 0 \to H^{0}(X, \Kx^{-1}L) \to H^{0}(X, \Kx ) \otimes H^{0}(X, L ) \xrightarrow{\mu} H^{0}(X, \Kx L ) \to \cdots \]
whence $\Ker ( \mu ) = H^{0}(X, \Kx^{-1}L)$, which is of dimension two by Riemann--Roch.
\par
We choose a basis $u$, $v$ for $H^{0}(X, F_{L}^{*} \otimes \Kx )$. This gives a map
\[ (u,v) \colon F_{L} \to \Kx \oplus \Kx. \]
This is of maximal rank: if $(u,v)$ factorised via a line bundle $M$ we would have $M = \Kx$, but $u$ and $v$ were chosen to be linearly independent. Now $(u,v)$ maps $F_L$ surjectively to a rank two subsheaf ${\mathcal G}$ of $\cans \oplus \cans$. Since $F_L$ is stable, $\mu({\mathcal G}) > 5/3$, and since ${\mathcal G}$ is a subsheaf of $\cans \oplus \cans$, we have $\mu({\mathcal G}) \leq 2$. Since $\mu({\mathcal G})$ is a half integer between $5/3$ and 2, it is equal to 2 and $(u,v)$ is surjective.
\par
Taking determinants, we find $\Ker(u,v) \cong \Kx^{-2}L$. Thus we have the sequence $0 \to \Kx^{-2}L \to F_{L} \to \Kx \oplus \Kx \to 0$. This completes the proof of Prop.\ \ref{defFLt}. \qed \\

We now show that, up to scalar, the above map $\Kx^{-2}L \to F_L$ is unique.

\begin{prop} The bundle $F_{L} \otimes \Kx^{2}L^{-1}$ has just one independent section. \label{onesympform} \end{prop}
\textbf{Proof}\\
Tensoring (\ref{defFL}) by $\Kx^{2}L^{-1}$ and taking global sections, we obtain
\[ 0 \to H^{0}(X, F_{L} \otimes \Kx^{2}L^{-1}) \to H^{1}(X, \Kx^{2}L^{-2}) \to H^{1}(X, \Kx^{2}L^{-1}) \otimes H^{0}(X, L)^{*} \]
and one checks that the second map is identified with the cup product map
\[ H^{1}(X, \Kx^{2}L^{-2}) \to \Homom \left( H^{0}(X, L) , H^{1}(X, \Kx^{2}L^{-1}) \right). \]
By Lemma \ref{multdualcup}, then, $H^{0}(X, F_{L} \otimes \Kx^{2}L^{-1})^*$ appears as the cokernel of the multiplication map $\mu \colon H^{0}(X, \Kx^{-1}L ) \otimes H^{0}(X, L) \to H^{0}(X, \Kx^{-1} L^{2} ) = \C^{7}$. 
To determine this cokernel, we consider two cases:
\paragraph{(a)} Suppose $|\Kx^{-1}L|$ is base point free. In this case there is an exact sequence of vector bundles $0 \to \Kx L^{-1} \to \Ox \otimes H^{0}(X, \Kx ^{-1} L ) \to \Kx ^{-1} L \to 0$, the surjection being the evaluation map. Tensoring by $L$ and taking global sections, we obtain
\begin{multline*} 0 \to H^{0}(X, \Kx) \to H^{0}(X , L) \otimes H^{0}(X, \Kx^{-1}L) \xrightarrow{\mu} \\ H^{0}(X , \Kx^{-1}L^{2}) \to H^{1}(X , \Kx ) \to 0 \end{multline*}
whence $\Coker (\mu) = H^{1}(X, \Kx)$, so is of dimension one.
\paragraph{(b)} If $\Kx^{-1}L = \Kx(x)$ for some $x \in X$ then $L = \Kx^{2}(x)$ and $H^{0}(X, \Kx^{-1}L ) = H^{0}(X, \Kx(x) )$. We have a diagram
\[ \xymatrix{ H^{0}(X, \Kx^{-1}L ) \otimes H^{0}(X, L) \ar[r]^-{\mu} \ar[d]^{\|} & H^{0}(X, \Kx^{-1}L^{2} ) \ar[d]^{\|} \\
H^{0}(X, \Kx(x) ) \otimes H^{0}(X, \Kx^{2}(x)) \ar[r]^-{\mu} & H^{0}(X, \Kx^{3}(2x) ) } \]
By the base point free pencil trick (Arbarello et al \cite{ACGH1985}, p.\ 126), the kernel of $\mu$ is isomorphic to $H^{0}(X, \Kx(x)) = \C^2$. Now $h^{0}(X, \Kx(x) ) \cdot h^{0}(X, \Kx^{2}(x)) = 2 \times 4 = 8$, so the image of $\mu$ is of dimension six. Thus $\Coker ( \mu )$ is again of dimension one.\\
\\
In either case, we obtain $h^{0}(X, F_{L} \otimes \Kx^{2} L^{-1}) = 1$. \qed \\

Now the symplectic form on $W$ induces a symplectic form on $W \otimes \Kx$ with values in $\Kx^2$. Furthermore, it is not hard to show (using for example Theorem \ref{symporthext}) that a subbundle $G \subset W$ is isotropic in $W$ if and only if $G \otimes \Kx$ is isotropic in $W \otimes \Kx$.

\begin{lemma} The degree one line bundle $\Kx^{-2}L$ is effective. \label{twentythree} \end{lemma}
\textbf{Proof}\\
Firstly, we claim that $h^{0}(X, W(x)) = 1$ for all $x \in X$. To see this, note that $H^{0}(X, W(x)) \cong H^{0}(X, W \otimes \Kx (-\iota x))$ and the latter space is identified with $\Ker \left( \Ox \otimes H^{0}(X, W \otimes \Kx )|_{\iota x} \to F_{L}|_{\iota x} \right)$. This map is surjective by Lemma \ref{Fpdesc}, so the kernel is of dimension one. This also shows that the generic value of $h^{0}(X, L \otimes W )$ as $L$ ranges over $\J$ is 1.
\par
By Thm.\ \ref{isotsubsh}, then, for each $x \in X$ which is not a Weierstrass point, the line subbundles $\Kx ( -x) = \Ox (\iota x)$ and $\Kx(-\iota x) = \Ox(x)$ define a rank two subsheaf of ${\cal W} \otimes \cans$, generating a rank two subbundle $G_x$ of $W \otimes \Kx$ on which the symplectic form inherited from $W \otimes \Kx$ vanishes identically\footnote{It is not hard to show, using for example Theorem \ref{symporthext}, that a subbundle $G \subset W$ is Lagrangian if and only if $G \otimes \Kx$ is Lagrangian with respect to the $\Kx^2$-valued symplectic form on $W \otimes \Kx$.}.
\par
We can say even more about the Lagrangian subbundles $G_x$. The line subbundles $\Ox(x) \subset W \otimes \Kx$ are contained in $F_L$, because an inclusion $\Ox(x) \hookrightarrow W \otimes \Kx$ is equivalent to a section $\Ox \to W \otimes \Kx$ vanishing at $x$, and all sections of $W \otimes \Kx$ are by definition $F_L$-valued. Therefore $G_x$ also belongs to $F_L$.
\par
We now give a more geometric way to realise the restriction of the $\Kx^2$-valued symplectic form on $W \otimes \Kx$ to $F_L$. Since $F_L$ has rank greater than two, the form does not restrict to zero (although it is degenerate).\\
\\
\textbf{Claim:} $h^{0}(X, \Homom(\bigwedge^{2}F_{L}, \Kx^{2})) = 1$. We have $\bigwedge^{2}F_{L} \cong \Homom(F_{L} , \det(F_{L}))$ since $F_L$ is of rank three. Therefore $\bigwedge^{2}F_{L}^{*} \otimes \Kx^{2} \cong F_{L} \otimes \Kx^{2} L^{-1}$, and this bundle has only one global section by Prop.\ \ref{onesympform}.\\
\par
Let $\alpha$ be a map $\Kx^{-2}L \to F_L$. Then
\begin{equation} w_{1} \wedge w_{2} \mapsto w_{1} \wedge w_{2} \wedge \alpha( \cdot ) \label{formonF}\end{equation}
defines a map $\bigwedge^{2} F_{L} \to \Homom ( \Kx^{-2}L , \det F_{L} ) \cong \Kx^2$. Clearly this is nonzero and, by the Claim, must be a scalar multiple of the restricted symplectic form.
\par
Therefore, the isotropy of $G_x$ in $W \otimes \Kx$ implies that the map
\[ \left( \Ox(\iota x) \oplus \Ox(x) \oplus \Kx^{-2}L \right) \to F_L \]
is not of maximal rank. In particular, for \emph{every} $x \in X$ apart from the Weierstrass points, the line subbundle $\Kx^{-2}L \subset F_L$ belongs to $G_x$.
\par
Since $F_L$ is stable, $G_x$ has degree two or three. If it is two then in fact $G_{x} = \Ox(\iota x) \oplus \Ox(x)$. Now $\Kx^{-2}L$ is also a degree one subbundle of $G_x$, so is equal to $\Ox(x)$ or $\Ox(\iota x)$.
\par
If $G_x$ has degree three then it is an elementary transformation
\[ 0 \to \strs(x) \oplus \strs(\iota x) \to {\cal G}_{x} \to \C_{a} \to 0 \]
for some $a \in X$. The determinant of $G_x$ is therefore $\Kx(a)$. Since $\Kx^{-2}L$ belongs to $G_x$, there is a diagram
\[ \xymatrix{ 0 \ar[r] & \Kx^{-2}L \ar[r] \ar[d]^{\|} & G_{x} \ar[r] \ar @{^{(}->}[d] & M \ar[r] \ar @{^{(}->}[d] & 0 \\
 0 \ar[r] & \Kx^{-2}L \ar[r] & F_{L} \ar[r] & \Kx \oplus \Kx \ar[r] & 0 } \]
where $M$ is a line bundle of degree two; clearly $M$ is isomorphic to $\Kx$. Therefore, since $\det(G_{x}) = \Kx(a)$, the bundle $\Kx^{-2}L$ is isomorphic to $\Ox(a)$.\\
In either case, $\Kx^{-2}L$ is effective. \qed \\

Note that $K^{-2}L \otimes \Kx^{-1} = \Ox(-\iota a)$. By Lemma \ref{defFLt}, the bundle $F_{L} \otimes \Kx^{-1}$ determines an extension class $e \in H^{1}(X, \Homom( \Ox \oplus \Ox , \Ox(- \iota a) ))$ which is of the form $(e_{1}, e_{2})$ for some $e_{1}, e_{2} \in H^{1}(X, \Homom( \Ox, \Ox(- \iota a)))$, a vector space of dimension two. Now we claim that $e_1$ and $e_2$ form a basis of this space. For, if $e_{2} = \mu e_1$ for some $\mu \in \C$ then the map $\Ox \to \Ox \oplus \Ox$ given by $\lambda \mapsto (-\mu \lambda, \lambda )$ would lift to $F_{L} \otimes \Kx^{-1}$, contradicting stability of this bundle.
\par
Now we consider a homomorphism $\Ox \to \Ox \oplus \Ox$ given by $\lambda \mapsto ( \alpha \lambda , \beta \lambda )$. The inverse image of the subbundle $(\alpha, \beta) \left( \Ox \right) \subset \Ox \oplus \Ox$ is an extension $E$ of $\Ox$ by $\Ox(- \iota a)$, and one sees, for example by inspecting transition functions, that the extension class of $E$ is
\[ \alpha e_{1} + \beta e_{2} \in H^{1}(X, \Homom( \Ox , \Ox(- \iota a) )). \]
Letting $(\alpha : \beta )$ vary in $\PP^1$, we find that all nontrivial extensions of this form are subbundles of $F_{L} \otimes \Kx^{-1}$ since $\{ e_{1}, e_{2} \}$ is a basis of the space
 $H^{1}(X, \Homom( \Ox , \Ox(- \iota a) ))$. They are also isotropic in $W$ because each one contains a pair of line subbundles of the form
\[ \Ox(-p), \quad \Ox(-\iota p) = \Kx^{-1}(p) \]
for some $p \in X$. Since $h^{0}(X, \Homom(\Ox(-x) , W)) = 1$ for all $x \in X$, by Thm.\ \ref{isotsubsh} these subbundles must generate an isotropic subbundle of rank two.

\begin{lemma} The point $a$ is a Weierstrass point of $X$. \end{lemma}
\textbf{Proof}\\
We have just seen that $W$ contains a pencil of isotropic subbundles which are nontrivial extensions $0 \to \Ox(- \iota a) \to E \to \Ox \to 0$. It will suffice to show that if $a$ is not a Weierstrass point, then any bundle containing such a pencil must have a theta divisor. To do this, we will show that for one of these $E$, the image of the map $\je \colon \J \dashrightarrow \pfed$ considered in section 5 is nondegenerate. By Lemma \ref{LtensW}, any symplectic bundle containing this $E$ as an isotropic subbundle must have a theta divisor.
\par
Suppose, then, that $\iota a \neq a$. We must show that one of the extensions $E$ above satisfies both of the genericity hypotheses of Lemma \ref{jnondeg}. Let $p \in X$ be a point such that $p \neq \iota p$, $p \neq a$ and $p \neq \iota a$ and consider the extension $E$ of $\Ox$ by $\Ox(-a)$ defined by the divisor
\[ a + p + \iota p \in |\Kx(a)| = \PP H^{1}(X, \Homom(\Ox, \Ox(-a))). \]
Then $E^*$ fits into a short exact sequence $0 \to M \to E^{*} \to \Ox(p) \to 0$ where $M := \Ox(a-p)$. Firstly, $M$ is not a point of order two in $\Pic^{0}(X)$. For, otherwise we would have $\Ox(2a) = \Ox(2p)$, whence either $a$ and $p$ are both Weierstrass points or $p=a$, both of which are excluded by our hypothesis.
\par
This means that there exists a unique pair of points $x,y$ such that $\Kx M^{2} = \Ox(x+y)$. Notice that $2a + \iota p \sim p + x + y$.
\par
It remains to check that $E$ satisfies the other genericity hypothesis of Lemma \ref{jnondeg}, so that $\je$ is defined at the points which would map to the point $\PP H^{0}(X, \Kx M^{2})$ in $\pfed$. We check as before that $\je(L)$ is this point if and only if $L^{-1} = M(-u)$ and $\Kx^{-1}L = M(-v)$, for some $u, v \in X$, and that the only solution to these equations (up to exchanging $L$ and $\Kx L^{-1}$) is
\[ L = M^{-1}(x) = M(\iota y) \quad \hbox{and} \quad \Kx L^{-1} = M^{-1}(y) = M(\iota x). \]
Now $\je$ is not defined at $M^{-1}(x)$ if and only if
\[ h^{0}(X, M^{-1}(x) \otimes E^{*}) \geq 2 \hbox{ and/or } h^{0}(X, M^{-1}(y) \otimes E^{*}) \geq 2, \]
equivalently, since $\chi(X, M^{-1}(x) \otimes E^{*}) = 1 = \chi(X, M^{-1}(y) \otimes E^{*})$,
\[ h^{1}(X, M^{-1}(x) \otimes E^{*}) \geq 1 \hbox{ and/or } h^{1}(X, M^{-1}(y) \otimes E^{*}) \geq 1. \]
By Serre duality, this becomes
\[ h^{0}(X, \Kx M(-x) \otimes E) \geq 1 \hbox{ and/or } h^{0}(X, \Kx M(-y) \otimes E) \geq 1, \]
that is, $\Kx^{-1} M^{-1}(x) = M^{-1}(- \iota x)$ or $\Kx^{-1} M^{-1}(y) = M^{-1}(-\iota y)$ belongs to $E$.
\par
Thus, in order that $\je$ be defined at $L$ (and thus also at $\Kx L^{-1}$), we require that the sets
\[ \{ M^{-1}(- \iota x), M^{-1}(- \iota y) \} \quad \hbox{and} \quad \{ \Ox(-a), \Ox(-p), \Ox(- \iota p) \} \]
be disjoint. Since $\Ox(- \iota p) = \Kx^{-1}(p)$, in fact we can replace the second one by $\{ \Ox(-a), \Ox(-p) \}$.
\par
Suppose $M^{-1}(- \iota x) = \Ox(-a)$. Then $\Ox(p - a - \iota x) = \Ox(-a)$, so $\Ox(p+a) = \Ox(\iota x + a)$. Thus $p = \iota x$. But then $\Ox(p+x+y) = \Kx(y)$. Since $2a + \iota p \sim p+x+y$, either $a$ is a Weierstrass point or $a = p$, contrary to hypothesis. The other possibilities can be excluded in a similar manner.
\par
Thus $\je$ is defined at the required points and so its image in $\pfed$ is nondegenerate. By Lemma \ref{LtensW}, for such an $E$, some symplectic extension of $E^*$ by $E$ has no theta divisor only if $a$ is a Weierstrass point. \qed \\

We have shown, therefore, that $F_{L} \otimes \Kx^{-1}$ contains all nontrivial extensions
\[ 0 \to \Ox(-a) \to E \to \Ox \to 0 \]
as isotropic subbundles. Since, by Lemma \ref{jnondeg}, there is at most one symplectic extension of each of $E^*$ by $E$ with no theta divisor, the bundle $W$ is isomorphic to the bundle $W_a$ constructed in section 6 (and we get another proof of the fact that $W_a$ contains all nontrivial extensions $E$ of the above form as Lagrangian subbundles, so the isomorphism class of $W_a$ depends only on $a$). Thus there are at most six possibilities for a bundle $W \in \Mt$ with no theta divisor.\\
\quad \par
It remains to show that the base locus of $|\Xi|$ is reduced.

\begin{lemma} Let $w, v \in X$ be distinct Weierstrass points. Then bundles $W_v$ and $W_w$ are mutually nonisomorphic. \end{lemma}
\textbf{Proof}\\
Let $E$ be a nontrivial extension of $\Ox$ by $\Ox(-w)$ defined by a divisor $p + \iota p + w \in |\Kx(w)|$ where $p \neq v$. Now $W_w$ is an extension of $E^*$ by $E$ by Lemma \ref{Windepe}. Consider now notrivial extensions of the form
\begin{equation} 0 \to \Ox(-v) \to F_{f} \to \Ox \to 0 \label{extOOv} \end{equation}
parametrised by $f:=\langle \delta(F) \rangle \in \PP H^{1} (X, \Homom(\Ox, \Ox(-v))) = \PP^1$. If $W_v$ were isomorphic\footnote{A priori, we should only ask for $S$-equivalence but, by Prop.\ \ref{bpstable}, the base locus of $|\Xi|$ consists of stable vector bundles so $S$-equivalence of $W_v$ and $W_w$ implies isomorphism.} to $W_w$ then we would have a map $F_{f} \hookrightarrow W_w$ for all such $F_f$, again by Lemma \ref{Windepe}. 
\par
There are three possibilities for the rank of the subsheaf ${\mathcal F}_{f} \cap {\mathcal E}$: these are zero, one and two. It is never two, because this would imply that some $F_f$ were isomorphic to $E$, which is impossible since $\Ox(-v) \subset F_f$ for all $f$ but $\Ox(-v) \not\subset E$. If it is one then it is not hard to check\footnote{See for example the proof of \cite{Hit2005b}, Prop.\ 18.} that the subsheaf ${\mathcal F}_{f} \cap {\mathcal E}$ corresponds to a line subbundle of degree $-1$ in $F_f$ and $E$, but $F_f$ and $E$ have such a subbundle in common for only finitely many $f$. Thus if $W_{v} \cong W_w$ then we have maps $F_{f} \to E^*$ of generic rank two for almost all $f \in \PP H^{1}(X, \Homom(\Ox, \Ox(-v)))$.
\par
By Riemann--Roch and since
\[ h^{1}(X, \Homom(\Ox(-v), E^{*})) = h^{0}(X, \Homom(\Ox(-v), E)) = 0, \]
there is one independent map $\gamma \colon \Ox(-v) \to E^*$. An extension $F_f$ admits a map to $E^*$ of generic rank two only if $\gamma$ factorises via $F_f$. By Narasimhan--Ramanan \cite{NR1969}, Lemma 3.2, this is equivalent to $\delta(F_{f})$ belonging to the kernel of the map
\[ H^{1}(X, \Homom( \Ox , \Ox (-v))) \to H^{1}(X, \Homom(\Ox, E^{*})) \]
induced by $\gamma$. Via Serre duality, this is dual to the map
\begin{equation}  H^{0}(X, \Kx \otimes E ) \to H^{0}(X, \Kx(v) ) \label{KEKv} \end{equation}
induced by $\gamma^{*} \colon E \to \Ox(v)$. We have rank two maps $F_{f} \to E^*$ for almost all $F_f$ if and only if (\ref{KEKv}) is zero. Now there is an exact sheaf sequence
\[ {\mathcal N} \to {\mathcal E} \xrightarrow{\gamma^*} \strs(v) \]
where ${\mathcal N}$ is invertible. The map in (\ref{KEKv}) is zero only if every map $\Kx^{-1} \to E$ factorises via the line subbundle $N^{\prime} \subset E$ generated by ${\mathcal N}$; this has degree at most $-1$. Now $\chi(X, \Homom(\Kx^{-1} , E)) = 1$ and
\[ h^{1}(X, \Homom(\Kx^{-1} , E)) = h^{0}(X, E^{*}) = 1, \]
so $h^{0}(X, \Homom(\Kx^{-1} , E)) = 2$. But $\Homom( \Kx^{-1} N^{\prime})$ is a line bundle of degree at most one, so there is at most one independent map $\Kx^{-1} \to N^{\prime}$. Therefore not every map $\Kx^{-1} \to E$ can factorise via $N^{\prime} \hookrightarrow E$.
\par
This means that there are no rank two maps $F_{f} \to E^*$ for most extensions $F_f$. Putting all this together, there cannot be maps $F_{f} \to W_w$ for all $F_f$. Hence $W_v$ cannot be isomorphic to $W_w$. \qed \\
\par
In summary, every symplectic bundle of rank four over $X$ with no theta divisor is of the form $W_w$ for some Weierstrass point $w$, and the bundles $W_w$ and $W_v$ are mutually nonisomorphic for $w \neq v$.\\
\\
This completes the proof of Theorem \ref{main}, our main result. \qed

\section{The link with Raynaud bundles}

In \cite{Ray1982}, Raynaud gives examples of semistable bundles without theta divisors over curves of arbitrary genus. Arnaud Beauville has shown that in the rank four / genus two case, some of these bundles admit symplectic or orthogonal structures. In this section, we show how the extensions in $\S$ 6 of the present article are related to Raynaud and Beauville's work.

\subsection{Raynaud's construction in genus 2}

This subsection is expository; the reference is Raynaud \cite{Ray1982}, sect.\ 3.\\ %We recall the details of Raynaud's construction in the genus two case.\\
\par
We write $J := J^0$ for brevity and identify $J$ with its dual Abelian variety $\hJ = \Pic^{0}(J)$ by means of the principal polarisation. We choose a symmetric divisor on $J$ defining the principal polarisation and, abusing notation, denote it $\Theta$. Consider the (ample) bundle $\Oj(2\Theta)$ on $J$. Via our identification $J \xrightarrow{\sim} \hJ$, the map $\phi_{2\Theta} \colon J \to \hJ$ is identified with the duplication map $2_{J} \colon J \to J$, which has degree $2^{2g} = 16$.\\
\par
Let $\mathbf{P}$ be the Poincar\'e bundle over $J \times J$ which is trivial over $\{ 0 \} \times J$ and $J \times \{ 0 \}$. We write $p$ and $q$ for the projections of $J \times J$ to the first and second factors.
\par
Recall that a sheaf $N$ over $J$ is \textsl{WIT} (``Weak Index Theorem'') of index $i$ if the sheaves $R^{j} q_{*} \left( p^{*} N \otimes \mathbf{P} \right)$ are zero for all $j \neq i$. Following Birkenhake--Lange \cite{BL2004}, p.\ 445, we define the \textsl{Fourier--Mukai transform} of such an $N$ as $R^{i} q_{*} \left( p^{*} N \otimes \mathbf{P} \right) =: F(N)$. Now $\strsj(-2\Theta)$ is WIT of index $g = 2$. In \cite{Muk1981}, Mukai proved:
\begin{enumerate}
\renewcommand{\labelenumi}{$\mathrm{(\roman{enumi})}$}
\item $F \left( \Oj(-2\Theta) \right)$ is a vector bundle $M$ of rank $(2 \Theta)^{g}/g = 4$ over $J$.
\item $F(M) \cong (-1_{J})^{*}\Oj(-2\Theta)$.
\item $2_{J}^{*}(M) \cong \Oj(2\Theta)^{\oplus 4}$.
\end{enumerate}
Let $U \subset J$ be an open set over which $\Oj(-2 \Theta)$ is trivial. Shrinking $U$ if necessary, we can assume that there is a section $s$ of $p^{*}M \otimes \mathbf{P}$ over $q^{-1}(U) = J \times U$ whose restriction to $\{ 0 \} \times U$ is nonzero. Consider the Abel--Jacobi map $\alpha_{c} \colon X \hookrightarrow J$ which sends the point $x \in X$ to $\Ox(x-c) \in J$. We denote the image curve $X_c$; clearly this passes through $0$. Then $\alpha_{c}^{*}M$ is a vector bundle $E_c \to X$ of rank four. By construction, the restriction of $s$ to $X_{c} \times U$ is nonzero, whence $h^{0}(X, E_{c} \otimes L) > 0$ for generic (and hence all) $L \in J$.
\par
To see that $E_c$ is semistable, we introduce the (possibly nonreduced) curve $Y := 2_{J}^{-1}(X_{c}) \subset J$. Write $f \colon Y \to X_c$ for the restriction of $2_J$ to $Y$, which is a degree 16 map of curves. Now for any vector bundle $V \to X_c$, we have
\begin{equation} \mu(f^{*}V) = \degree(f) \cdot \mu(V) = 16 \: \mu(V). \label{slopes} \end{equation}
By (iii), the pullback of $E_c$ to $Y$ is a direct sum of line bundles of the same degree, so is semistable. From this and (\ref{slopes}) we check that $E_c$ is also semistable.
\par
Let us find the slope of $E_c$. Firstly, we have
\[ \degree ( E_{c} ) =\frac{\degree ( (2_{J}^{*}M)|_{Y} )}{\degree (f)}. \]
On the other hand, since $(2_{J})^{*}M \cong \Oj(2\Theta)^{\oplus 4}$, we have
\begin{align*} \degree ( (2_{J}^{*}M)|_{Y} ) &= 4 \: \degree ( \Oj(2\Theta)|_{Y} )\\
&= \degree ( \Oj(2\Theta)|_{Y}^{\otimes 4} )\\
&= \degree ( ( (2_{J})^{*}\Oj(2\Theta))|_{Y} ) \end{align*}
since $(n_{J})^{*}\Oj(2\Theta) = \Oj(2n^{2}\Theta)$. This is in turn equal to
\[ \degree (2_{J}) \cdot (\degree (\Oj(2\Theta)|_{X_c}) = 16 ( 2\Theta \cdot X_{c} ) = 16 \cdot 2g = 64, \]
so $f^{*}E_c$ has slope 16. By (\ref{slopes}), the slope of $E_c$ is 1.\\
\\
\emph{Thus, the tensor product of $E_c$ by any line bundle of degree $-1$ is a semistable vector bundle of degree zero and rank four with no theta divisor.}\\
\\
\textbf{Remark:} In fact this is a stable bundle, since, by Raynaud \cite{Ray1982}, Cor.\ 1.7.4, all semistable bundles of degree zero and rank at most three over a curve of genus two have a theta divisor.

\subsection{Symplectic and orthogonal structures}

Here we show that for certain $c \in X$, the bundle $E_c$ admits a symplectic structure.

\paragraph{Note:} This construction is due to Arnaud Beauville, to whom I express my thanks for allowing me to present his results here.\\
\par
We begin by describing three actions of the Heisenberg group $\Htt$. Set-theoretically, this group consists of pairs $(\phi, \eta)$ where $\eta$ is a point of order two in $J$ and $\phi$ an automorphism of the variety $\Oj(2\Theta)$ covering $t_\eta$; in other words, such that there is a commutative diagram of varieties
\[ \xymatrix{ \Oj(2\Theta) \ar[r]^{\phi} \ar[d] & \Oj(2\Theta) \ar[d] \\ J \ar[r]^{t_{\eta}} & J. } \]
Write $J[2] := \Ker (2_{J})$. The group $\Htt$ is a central extension
\[ 1 \to \C^{*} \to \Htt \to J[2] \to 0 \]
by for example Birkenhake--Lange \cite{BL2004}, Prop.\ 6.1.1 and since $\phi_{n\Theta}$ is multiplication by $n$ on $J$.
\par
(i) The first action of $\Htt$ is on $H^{0}(J, 2\Theta)$ and is given in \cite{BL2004}, $\S$ 6.4. For a global section $s$ of $\Oj(2\Theta)$, we define $(\phi , \eta) \cdot s = \phi \circ s \circ t_{-\eta}$, which is equal to $\phi \circ s \circ t_{\eta}$ since $-\eta = \eta$. This gives a representation of $\Htt$ in $\mathrm{GL}( H^{0}(J, \Oj(2\Theta)))$.\\
\par
(ii) There is a natural action of $\Htt$ on the total space of $\Oj(-2\Theta)$, given by $(\phi , \eta ) \cdot k = (^{t}\phi^{-1})(k)$.\\
\par
(iii) The last action will define a linearisation of $J[2]$ on $\Oj(4\Theta)$, which will depend on the choice of a theta characteristic $\kappa$. To each $\kappa$, we associate a character of weight $2$ of $\Htt$, denoted $\chi_{\kappa}$ (see Beauville \cite{Bea1991}, Lemme A.4). Given $(\phi , \eta) \in \Htt$, the map
\begin{equation} (\phi \otimes \phi) \cdot \chi_{\kappa}(\phi, \eta)^{-1} \label{action4Theta} \end{equation}
is an automorphism of $\Oj(4\Theta)$ covering $t_{\eta}$. The subgroup $\C^{*} \subset \Htt$ acts trivially, so this action factorises via $\Htt \to J[2]$.\\
\par
Fix a theta characteristic $\kappa$, and consider the product of these three actions on the bundle
\[ Q := \Homom \left( \strsj(2\Theta) \otimes H^{0}(J, 2\Theta)^{*} , \left(\strsj(2\Theta) \otimes H^{0}(J, 2\Theta)^{*} \right)^{*} \otimes \strsj(4\Theta) \right). \]
We notice firstly that the actions of $\C^*$ on each copy of $\Oj(-2\Theta) \otimes H^{0}(J, 2\Theta)$ cancel. Since the subgroup $\C^*$ acts trivially on $\Oj(4\Theta)$, the action of $\Htt$ on $Q$ factorises via $\Htt \to J[2]$.
\par
Now by Schneider \cite{Sch2005}, Prop.\ 2.1, there is an isomorphism of $J[2]$-bundles between $\Oj(2\Theta) \otimes H^{0}(J, \Oj(2\Theta))^*$ and $2_{J}^{*}M$, so $\Oj(2\Theta) \otimes H^{0}(J, \Oj(2\Theta))^*$ descends to $M$. Furthermore, we have

\begin{prop}
The quotient of $\Oj(4\Theta)$ by the linearisation of the action of $J[2]$ corresponding to the theta characteristic $\kappa$ is $\Oj(\Theta_{\kappa})$. \end{prop}
\textbf{Proof}\\
The line bundles over $J$ which pull back to $\Oj(4\Theta)$ are exactly the $\Oj(\Theta_{\kappa})$ by for example Birkenhake--Lange, p.\ 34, and since $\Ker( 2_{J}^{*} ) \cong J[2]$. Therefore, the quotient, if it exists, is among the $\Oj(\Theta_{\kappa})$.
\par
For each $(\phi, \eta) \in \Htt$, the map $\phi^{\otimes 2}$ is an automorphism of $\Oj(4\Theta)$ covering $t_{\eta}$. Thus we get a homomorphism $\Htt \to {\mathcal G}(4\Theta)$ and an action of $\Htt$ on $H^{0}(J, \Oj(4\Theta) )$.
\par
Now each $\Oj(\Theta_{\kappa})$ has a unique global section up to scalar since the polarisation is principal. We write $\vartheta_{\kappa}$ for this section. By Beauville \cite{Bea1991}, Prop.\ A.8, the pullbacks of the $\vartheta_{\kappa}$ by $2_J$ form a basis for the space $H^{0}(J, \Oj(4\Theta) )$. Furthermore, for each $\kappa$, the section $2_{J}^{*} \vartheta_{\kappa}$ is an eigenvector for the action of $\Htt$ via $\chi_{\kappa}$; precisely,
\[ (\phi, \eta) \cdot (2_{J}^{*} \vartheta_{\kappa}) =  \chi_{\kappa}(\phi, \eta) (2_{J}^{*} \vartheta_{\kappa}). \]
\par
Now we will fix a theta characteristic $\kappa$ and twist this action of $\Htt$ on $\Oj(4\Theta)$ by $\chi_{\kappa}^{-1}$, so we get the linearisation of $J[2]$ on $\Oj(4\Theta)$ in (iii) above. Since the underlying action of $J[2]$ on $J$ is free, by Kempf's Descent Lemma $\Oj(4\Theta)$ can also be descended to $J \cong J / J[2]$ by this linearisation. The induced action of $J[2]$ on $H^{0}(J, 4\Theta)$ is given by
\[ (\phi, \eta) \cdot s = \chi_{\kappa}(\phi, \eta) ^{-1} \left( \phi^{\otimes 2} \circ s \circ t_{\eta} \right). \]

Putting all this together, we see that $2_{J}^{*} \vartheta_{\kappa}$ is equivariant for the linearisation corresponding to $\kappa$. Thus it descends to the quotient by this linearisation. It is not hard to see that the descended section vanishes exactly along $\Theta_{\kappa}$ and that the quotient is $\Oj(\Theta_{\kappa})$. \qed \\
\par
Thus, for a theta characteristic $\kappa$, the bundle $Q$ descends via the linearisation corresponding to $\kappa$ to $\Homom(M, M^{*} \otimes \Oj(\Theta_{\kappa}))$. We now describe a section of $Q$ which is equivariant for this linearisation. By Beauville \cite{Bea1991}, Prop.\ A.5, there is a basis of
\[ \Homom_{\C} ( H^{0}( J, 2\Theta )^{*} , H^{0}(J, 2 \Theta) ) = H^{0}( J, 2\Theta ) \otimes H^{0}( J, 2\Theta ) \]
indexed by the theta characteristics of $X$. An element of this basis is an eigenvector $\xi_{\kappa}$ with respect to the character $\chi_{\kappa}$ of weight 2:
\begin{equation} (\phi , \eta ) \cdot \xi_{\kappa} = \chi_{\kappa}(\phi , \eta) \xi_{\kappa}. \label{A} \end{equation}
Clearly the bundle $Q$ is isomorphic to $\Oj \otimes H^{0}(J, 2\Theta)^{\otimes 2}$, so each $\xi_{\kappa}$ defines a section $1 \otimes \xi_{\kappa}$ of $Q$. Let $j \in J$; then one checks that an automorphism $(\phi, \eta) \in \Htt$ sends the element $1_{j} \otimes \xi_{\kappa} \in \Oj|_{j} \otimes H^{0}( J, 2\Theta )^{\otimes 2} = Q|_j$ to
\[ \left( \chi_{\kappa}(\phi, \eta)^{-1} \cdot 1_{j + \eta} \right) \otimes \left( \chi_{\kappa}(\phi, \eta) \xi_{\kappa} \right) = 1_{j + \eta} \otimes \xi_{\kappa} \in Q|_{j + \eta}. \]
Thus $1 \otimes \xi_{\kappa}$ is an equivariant section for the linearisation of $J[2]$ on $Q$ corresponding to $\kappa$. Hence it descends to $\Homom(M, M^{*} \otimes \Oj(\Theta_{\kappa}))$. Moreover, $\xi_{\kappa}$ defines an isomorphism $H^{0}(J, 2\Theta)^{*} \to H^{0}(J, 2\Theta)$ by \cite{Bea1991}, Remarque A.6, so $1 \otimes \xi_{\kappa}$ defines an isomorphism $M \xrightarrow{\sim} M^{*} \otimes O_{J}(\Theta_{\kappa})$. Furthermore, by \cite{Bea1991}, Prop.\ A.5, this isomorphism is symmetric or antisymmetric according to the parity of $\kappa$. In this way we get get a symplectic (resp., orthogonal) structure on $M$ for each odd (resp., even) theta characteristic.\\
\par
Now let $w \in X$ be a Weierstrass point and consider the Abel--Jacobi map $\alpha_{w} \colon X \to J$. We write $\tau := O_{X}(w)$. Denote $E_w$ the bundle $\alpha_{w}^{*}M \to X$. By the above argument, there exists an antisymmetric isomorphism
\[ 1 \otimes \xi_{\tau} \colon E_{w} \to E_{w}^{*} \otimes ( \Oj(\Theta_{\tau})|_{X_w}). \]
Now $K_J$ is trivial and $X_{w} = \Theta_{\tau}$ as divisors. Therefore, by the adjunction formula,
\[ \Oj(\Theta_{\tau})|_{X_w} \cong K_{J} \otimes \Oj(X_{w})|_{X_w} \cong \Kx, \]
so $E_w$ carries a $\Kx$-valued symplectic form. A straightforward calculation shows that $\alpha_w$ is the only Abel--Jacobi map which induces in this way a $\Kx$-valued symplectic form on $\alpha_{w}^{*}M$. We set $V_{w} = E_{w} \otimes \kappa^{-1}$ for any theta characteristic $\kappa$; we will see shortly that the isomorphism class of $V_w$ does not depend on this choice of $\kappa$. By the preceding discussion, we have an antisymmetric isomorphism $V_{w} \xrightarrow{\sim} V_{w}^*$.\\
\\
This construction has an interesting corollary:
\begin{prop} The bundle $V_w$ is invariant under tensoring by all line bundles of order two in $\Pic^{0}(X)$. \end{prop}
\textbf{Proof}\\
Composing the isomorphisms $M^{*} \otimes \Oj(\Theta_{\kappa}) \xrightarrow{\sim} M$ and $M \xrightarrow{\sim} M^{*} \otimes \Oj(\Theta_{\kappa^{\prime}})$ for each pair $\kappa$, $\kappa^{\prime}$ of distinct theta characteristics, we see that $M \cong M \otimes \eta$ for each two-torsion point $\eta$ in $J$. By Birkenhake--Lange \cite{BL2004}, Lemma 11.3.1, the map $\alpha_{w}^{*} \colon \Pic^{0}(J) \to \Pic^{0}(X)$ is an isomorphism, so $E_w$ and $V_w$ are also invariant under tensoring by line bundles of order two. \qed \\
\\
This proposition shows in particular that $E_{w} \otimes \kappa^{-1}$ and $E_{w} \otimes (\kappa^{\prime})^{-1}$ are both isomorphic to $V_w$, as we needed.\\
\\
\emph{In summary, we obtain a stable rank four symplectic bundle without a theta divisor for each odd theta characteristic of $X$.}

\subsection{The link with extensions}

Here we make explicit the link between the bundles of Raynaud and Beauville just described and those constructed in section 6. The result, philosophically the only one possible, is that for each Weierstrass point $w$, the bundles $V_w$ and $W_w$ are isomorphic. This will follow easily from Lemma \ref{xixw}.% 
\begin{lemma} For any $x \in X$, the fibres of the subbundles $\Ox(-x)$ and $\Ox(-\iota x)$ of $V_w$ coincide at the point $w$. \label{xixw} \end{lemma}
\textbf{Proof}\\
Since $F(M) = R^{0}q_{*} (p^{*}M \otimes \mathbf{P})$ is isomorphic to $\Oj(-2\Theta)$, by adjunction (Hartshorne \cite{Har1977}, p.\ 110) we have a map $\Psi \colon q^{*}\Oj(-2\Theta) \to p^{*}M \otimes \mathbf{P}$ over $J \times J$. This has the property that for each $(L_{1}, L_{2}) \in J \times J$, the image of $\Psi|_{(L_{1}, L_{2})}$ in
\[ ( p^{*}M \otimes \mathbf{P} ) |_{(L_{1}, L_{2})} = (M \otimes L_{2})|_{L_1} \]
is identified with that of the (unique) global section of $M \otimes L_{2}$ at $L_1$.
\par
We restrict $\Psi$ to the subvariety $\{ 0 \} \times J$ and denote $\Psi_0$ the induced map of vector bundles over $J$:
\[ \xymatrix @C=0.8in{ \left( q^{*}\Oj(-2\Theta) \right) |_{ \{ 0 \} \times J} \ar[r]^-{\Psi|_{ \{ 0 \} \times J}} \ar[d]^{\|} & \left( p^{*}M \otimes \mathbf{P} \right) |_{ \{ 0 \} \times J} \ar[d]^{\|} \\
 \Oj(-2\Theta) \ar[r]^-{\Psi_0} & M|_{0} \otimes \left( \mathbf{P}|_{ \{ 0 \} \times J} \right). } \]
Now recall that $\mathbf{P}$ is trivial over $\{ 0 \} \times J$, so $\Psi_0$ defines a section of $\Oj(2\Theta)^{\oplus 4}$. Since every section of $\Oj(2\Theta)$ is $-1_J$-invariant, we have a commutative diagram
\[ \xymatrix @C=0.8in{ \Oj(-2\Theta) \ar[r]^-{\Psi_0} \ar[d]^{\tmo} & M|_{0} \otimes \left( \mathbf{P}|_{ \{ 0 \} \times J} \right) \ar[d]^{\Iden_{M|_0} \otimes l} \\
\Oj(-2\Theta) \ar[r]^-{\Psi_0} & M|_{0} \otimes \left( \mathbf{P}|_{ \{ 0 \} \times J} \right) } \]
where $\tmo$ and $l$ are suitable linearisations of $\Oj(-2\Theta)$ and $\Oj$ respectively. This shows that the image of $\Psi_0$ at $L$ coincides with that of $L^{-1}$ under the identification of projective spaces $\PP ( M|_{0} \otimes L ) = \PP ( M|_{0} \otimes L^{-1} )$.
\par
In particular the fibres of $\Ox(x - w)$ and $\Ox(x - w)^{-1} = \Ox(\iota x - w)$ coincide in the fibre of the restricted bundle $E_w$, for each $x \in X$. This implies that the fibres of the subbundles $\Ox(-x)$ and $\Ox(- \iota x)$ of $E_{w}(-w) = V_w$ coincide in the fibre over $\alpha_{w}^{-1}(0) = w$, as required. \qed

\begin{thm} The bundle $V_w$ associated to a Weierstrass point $w$ is isomorphic to the extension $W_w$ constructed in section 6. \end{thm}
\textbf{Proof}\\
By Lemma \ref{xixw}, the fibres of the subbundles $\Ox(-x)$ and $\Ox(-\iota x)$ of $V_w$ coincide at the point $w$. Since $V_w$ is stable, they do not coincide anywhere else since then they would generate a rank two subbundle of degree at least zero. Thus $V_w$ contains an subbundle of rank two which is an elementary transformation
\[ 0 \to \strs(-x) \oplus \strs(- \iota x) \to {\mathcal E} \to \C_{w} \to 0 \]
Since $E$ has degree $-1$ and belongs to $V_w$, it is stable. Then it is not hard to see that it is one of the bundles $E_e$ that we considered in section 6. Moreover, it contains the subbundles $L^{-1}$, $\Kx^{-1} L$ for some $L = \Ox(x) \in \J$. Since we saw that $h^{0}(X, V_{w}(x)) =1$ for all $x \in X$, Theorem \ref{isotsubsh} shows that $E$ is isotropic in $V_w$. Thus $V_w$ is a symplectic extension of $E^*$ by $E$. Since by Lemma \ref{jnondeg} there is only one such extension with no theta divisor, $V_w$ must be isomorphic to $W_w$. \qed \\

This establishes the link between the two constructions of the base points of $|\Xi|$.

\section{Applications and future work}

We give one immediate application of the finiteness of $\Bs|\Xi|$.

\begin{cor} The map $D \colon \Mt \dashrightarrow \ftp$ is surjective. \end{cor}
\textbf{Proof}\\
Let $G$ be an even $4\Theta$ divisor. Choose nine hyperplanes $\Pi_{1}, \ldots, \Pi_{9}$ in $\ftp$ whose intersection consists of the point $G$. The inverse images $D^{-1}(\Pi_{1}), \ldots,$ $D^{-1}(\Pi_{9})$ intersect in a subset $S$ of dimension at least 1 in $\Mt$ since this variety is of dimension 10. Now $S$ contains the base locus of $|\Xi|$. By Theorem \ref{main}, this is finite, so the map $D$ is defined at most points of $S$. By construction, the image of any of these points is the divisor $G$. \qed \\

We finish with a logically independent remark. It is not hard to see that the theta divisor of a \emph{strictly semistable} bundle of degree zero, if it exists, must be reducible. The converse, however, is not true: there exist \emph{stable} bundles with reducible theta divisors, as studied by Beauville \cite{Bea2003}. In the case at hand, we can say:

\begin{prop} Suppose $X$ has genus two. Then there exist \emph{stable} symplectic bundles in $\Mt$ with reducible theta divisors. \end{prop}
\textbf{Proof}\\
Let $V_1$ and $V_2$ be mutually nonisomorphic stable bundles of rank two and trivial determinant over $X$ with (irreducible) theta divisors $D_{1}, D_{2} \in |2\Theta|$. Narasimhan and Ramanan showed in \cite{NR1969} that the theta map
\[ D \colon \mathcal{SU}_{X}(2, \Ox) \to |2\Theta| = |2\Theta|_+ \]
is an isomorphism. Hence $D_{1} \neq D_2$. Let $W$ be a strictly semistable bundle with $D(W) = D_{1} + D_2$. With the list on p.\ 25 of \cite{Hit2005}, we consider all the possibilities for $W$ and, using \cite{NR1969}, we find that the only one which can have a theta divisor of the form $D_{1} + D_2$ is a direct sum of two stable bundles $W_1$ and $W_2$ of rank 2 and trivial determinant. We have
\[ D(W_{1} \oplus W_{2}) = D(W_{1}) + D(W_{2}) = D_{1} + D_2 \]
so, since $\J$ is normal, $\{ D(W_{1}), D(W_{2}) \} = \{ D_{1}, D_{2} \}$.
\par
But $\mathcal{SU}_{X}(2, \Ox)$ is isomorphic to $|2\Theta|$, so $W_{1} \oplus W_{2} \cong V_{1} \oplus V_2$. Thus the fibre of $D$ over $D( V_{1} \oplus V_{2} ) = D_{1} + D_2$ contains at most one strictly semistable bundle. Since this fibre is of dimension at least one, there exist stable bundles with theta divisor $D_{1} + D_2$. \qed

\subsection*{Future work}

By Criterion \ref{symporthext}, one can construct bundles with orthogonal structure as extensions $0 \to E \to V \to E^{*} \to 0$ with classes in $H^{1}(X, \bigwedge^{2}E)$. This space again carries an action of $\iota$, so one might try to find orthogonal bundles without theta divisors with a construction analogous to that in section 6. One might also expect this construction to generalise to bundles of higher rank over hyperelliptic curves of higher genus.
\par
Other natural things to look for include a description of the fibres of $D \colon \Mt \dashrightarrow \ftp$ and an explicit construction of the stable bundles with reducible theta divisors just mentioned. These questions will be studied in the future.\\

\paragraph{Acknowledgements:} It is a pleasure to thank my doctoral supervisor Christian Pauly for all his time and interest, and for much useful discussion and advice. His influence on this work is considerable.
\par
My sincere thanks to Arnaud Beauville for his help and interest. In addition to answering many questions, he showed me in particular Prop.\ \ref{count} and the construction in $\S$ 9.2.
\par
I am grateful to Andr\'e Hirschowitz and Klaus Hulek for many helpful discussions. Thanks to my supervisors John Bolton, Wilhelm Klingenberg and Bill Oxbury at Durham for all their mathematical and practical help during my doctoral studies.
\par
I acknowledge gratefully the financial and practical support of the Universities of Durham, Nice--Sophia-Antipolis and Hannover, the Engineering and Physical Sciences Research Council and the Deutsche Forschungsgemeinschaft Schwerpunktprogramm, ``Globale Methoden in der Komplexen Geometrie''.

\appendix

\section{Proof of Proposition \ref{Fpdesc}}

\textbf{Proposition \ref{Fpdesc}} \textsl{The evaluation map $\ev \colon \Ox \otimes H^{0}(X, W \otimes \Kx) \to W \otimes \Kx$ is everywhere surjective to a stable subbundle of rank three and degree five in $W \otimes \Kx$.}\\
\textbf{Proof}\\
The image of $\ev$ is a subsheaf ${\mathcal F}$ of ${\mathcal W} \otimes \cans$. We find first the rank of the corresponding vector bundle $F$, which is not necessarily a subbundle of $W \otimes \Kx$. Firstly, since we have seen that for any $x \in X$ there is a section of $W \otimes \Kx$ vanishing at $x$, the rank of $F$ cannot be four.
\par
We have $H^{0}(X, F) \cong H^{0}(X, W \otimes \Kx)$, so $h^{0}(X, F) = 4$. Since $W \otimes \Kx$ is stable, $\mu(F) < 2$. But no line bundle of degree one on $X$ has four independent sections, so $F$ is not of rank one. Suppose the rank is two. Then $\degree (F) \leq 3$. We now recall that it contains a subbundle $\Ox(x)$ for every $x \in X$. This is because it has sections vanishing at each point of the curve, but has none vanishing at more than one point because such a section would generate a line subbundle of degree at least two in $W$. Thus we have a short exact sequence $0 \to \Ox(x) \to F \to M \to 0$ for each $x \in X$, where $M$ is a line bundle of degree at most two (depending on $x$). Then
\[ h^{0}(X, F) \leq h^{0}(X, \Ox(x)) + h^{0}(X, M) \leq 1 + 2 = 3, \]
a contradiction.
\par
Thus $F$ has rank three. By stability of $W \otimes \Kx$, it has degree at most five. Suppose it is four or less. Let $x \in X$ be such that $\iota x \neq x$. As above, we have a short exact sequence $0 \to \Ox(x) \to F \xrightarrow{q} G \to 0$ where $G$ is of rank two and degree at most three.
\par
By hypothesis, $h^{0}(X,G) \geq 3$. An upper bound for the degree $d$ of a line subbundle of $G$ can be calculated as follows: the inverse image in $F$ of such a subbundle is a bundle of rank two and degree $d+1$. By stability of $W$, we have $\frac{d+1}{2} < 2$, so $d \leq 2$.
\par
Now let $t$ be a section of $F$ vanishing at $y \neq x$. The composed map $q \circ t$ generates a subbundle $M$ of degree one or two in $G$, so we have an exact sequence
\[ 0 \to M \to G \to N \to 0 \]
where $N$ is a line bundle of degree at most $3 - \degree(M)$. Now
\[ 3 \leq h^{0}(X, G) \leq h^{0}(X, M ) + h^{0}(X, N ). \]
This shows that one of $M$ and $N$ is the canonical bundle and the other is effective (in particular, $G$ has degree three). Suppose $N = \Kx$; then $M = \Ox(y)$. But we can do this for any $y$ (with the same $x$), which implies that $\det(G) = \Kx(y)$ for all $y \in X$, which is clearly absurd. Hence $M = \Kx$ and $N = \Ox(p)$ for some $p \in X$. We write $p = p_x$ since $p$ depends on $x$. Now $\det(F) = \Kx(x+p_{x}) = \Kx(x^{\prime}+p_{x^{\prime}})$ for generic $x^{\prime} \in X$. Therefore $\Ox(x+p_{x}) = \Ox(x^{\prime}+p_{x^{\prime}})$ for all $x, x^{\prime} \in X$, so $p_{x} = \iota x$ and $G$ is an extension $0 \to \Kx \to G \to \Ox (\iota x) \to 0$.
\par
Now by hypothesis, the unique section of $\Ox(\iota x)$ lifts to $G$, so the extension class $\delta(G)$ belongs to the kernel of the induced map
\[ H^{1}(X , \Homom(\Ox(\iota x) , \Kx)) \to H^{1}(X , \Homom(\Ox , \Kx)), \]
which is identified with $H^{0}(X , \Ox(\iota x))^{*} \to H^{0}(X , \Ox)^*$ by Serre duality. But this is an isomorphism, so $G = \Kx \oplus \Ox(\iota x)$.
\par
This means that $F$ is an extension $0 \to \Ox(x) \to F \to \Kx \oplus \Ox(\iota x) \to 0$ of class $\delta(F) \in H^{1}(X, \Homom( \Kx \oplus \Ox(\iota x) , \Ox(x) ) )$, that is,
\[ H^{1}(X, \Homom( \Kx , \Ox(x) ) ) \oplus H^{1}(X, \Homom( \Ox(\iota x) , \Ox(x) ) ) = \C^{3}. \]
Now $h^{0}(X, \Kx \oplus \Ox(\iota x)) = 3$, so all sections of $\Kx \oplus \Ox(\iota x)$ lift to $F$ by hypothesis. This means that the cup product map
\[ \cdot \cup \delta(F) \colon H^{0}(X, \Kx \oplus \Ox(\iota x) ) \to H^{1}( X, \Ox(x) ) \]
is zero. But by Lemma \ref{multdualcup}, the linear map
\begin{multline*} \cup \colon H^{1}(X, \Homom( \Kx \oplus \Ox(\iota x) , \Ox(x) ) ) \to \\ \Homom \left( H^{0}(X, \Kx \oplus \Ox(\iota x) ) , H^{1}( X, \Ox(x) ) \right) \end{multline*}
is dual to the multiplication map
\[  H^{0}(X, \Kx \oplus \Ox(\iota x) ) \otimes H^{0}( X, \Ox(\iota x) ) \to H^{0}(X, \Kx(\iota x) \oplus \Ox(2 \iota x )) \]
which is an isomorphism. Thus $\cup$ is also an isomorphism, so $\delta(F) = 0$. But then $\Kx$ is a subbundle of $F$, which is not possible since $W \otimes \Kx$ is stable.\\
\\
The only possibility, then, is that $F$ has degree five. To see that it is stable, note that by stability of $W \otimes \Kx$, any line subbundle of $F$ has slope at most 1, and any rank two subbundle has slope at most 3/2. Both of these are less than $5/3 = \mu(F)$, so $F$ is stable.\\
\par
Lastly, let $F^{\prime}$ be the subbundle of $W \otimes \Kx$ generated by the image of $\ev$. This has degree at least five since ${\mathcal F}$ has degree five. By stability of $W$, it is at most five. Hence $F^{\prime} = F$ and $\ev$ is everywhere of rank three. This completes the proof. \qed

\noindent Institut f\"ur Algebraische Geometrie,\\
Leibniz Universit\"at Hannover,\\
Welfengarten 1,\\
30167 Hannover,\\
Germany.\\
Email: \texttt{hitching@math.uni-hannover.de}

\end{document}